\documentclass[11pt]{article}
\usepackage{amsmath,amssymb,amsthm}
\newtheorem{thm}{Theorem}[section]

 \newtheorem{lem}{Lemma}[section]
 \newtheorem{prop}{Proposition}[section]
 \newtheorem{defn}{Definition}[section]
\theoremstyle{remark}
\newtheorem{rem}{Remark}[section]

 \allowdisplaybreaks
\numberwithin{equation}{section}

\title{\textbf{Global exponential stability of classical \\ solutions to the hydrodynamic \\
model for semiconductors}}
\author{Daoyuan Fang\thanks {E-mail: dyf@zju.edu.cn}\ \ \ \ Jiang Xu\thanks {E-mail: jiangxu\underline{ }79@yahoo.com.cn}\ \ \ \ Ting Zhang\thanks{E-mail: zhangting79@hotmail.com }\\
\small{\textit{Department of Mathematics, Zhejiang University},}
\\\small{\textit{Hangzhou 310027, P.R.China}}}
\date{}
\begin{document}
\maketitle{} \begin{abstract}
 \small{ In this paper, the global well-posedness and stability
of classical solutions to the multidimensional hydrodynamic model
for semiconductors on the framework of Besov space are considered.
We \textit{weaken} the regularity requirement of the initial data,
and improve some known results in Sobolev space. The local
existence of
  classical solutions to the Cauchy problem is obtained by the regularized means and compactness argument.
  Using the high- and low- frequency
decomposition method, we prove the global exponential stability of
classical solutions (close to equilibrium). Furthermore, it is
also shown that the vorticity decays to zero exponentially in the
$2$D and $3$D space. The main analytic tools are  the
Littlewood-Paley decomposition and Bony's para-product
 formula.}.\end{abstract}

\hspace{-0.5cm}\textbf{Keywords:} \small{Hydrodynamic; exponential
stability; classical solutions;
spectral localization.}\\

\hspace{-0.5cm}\textbf{AMS subject classification:} \small{35L65,
76X05, 35M10}

\section{Introduction and Main Results}
In the modern semiconductors industry, the numerical simulation of
device modeling has become very important. With the fast development
of miniaturization devices, the traditional drift-diffusion model is
no more valid, especially in submicron devices or in the occurrence
of high field phenomena. Some kinetic models such as Boltzmann
equation which  describe the evolution of the distribution function
$f(t,x,\textbf{v})$ of charged particles in the phase space are more
accurate, but they need much computing power in practical
application. Therefore, the hydrodynamic model which represents a
reasonable compromise between the physical accuracy and the
reduction of computational cost has recently received increasing
attention in statistical physics and applied mathematics. By
applying a moment method and appropriate closure conditions, it can
been derived from the semiclassical Boltzmann equation coupled with
the electric potential through a Poisson equation (see, e.g.,
Refs.~\cite{J}). In this paper, we are interested in the simplified
hydrodynamic model where the energy equation is replaced by a
pressure-density relation, which was first analyzed by Degond and
Markowich \cite{De}. The main objective is to study the global
well-posedness and stability of classical solutions to the
(unipolar) hydrodynamic model on the framework of Besov space. After
an appropriate scaling, it can be written as
$$\begin{cases}
n_{t}+\nabla\cdot(n\textbf{u})=0\cr
 \textbf{u}_{t}+(\textbf{u}\cdot\nabla)\textbf{u}+\frac{1}{n}\nabla p(n)=\nabla\Phi-\frac{\textbf{u}}{\tau} \cr
 \Delta \Phi=n-\bar{n},\ \ \ \Phi\rightarrow 0 \ \ \mbox{as} \ \ |x|\rightarrow +\infty
 \end{cases} \eqno(1.1)$$
 for $ (t,x)\in[0,+\infty)\times\mathbb{R}^{N}, \ N\geqslant2$,
 where $n,\textbf{u}=(u^1,u^2,\cdot\cdot\cdot,u^{N})^{\top}$($^\top$ represents transpose) and $\Phi$ denote the electron density, the electron velocity
 and the electrostatic potential respectively. $\nabla$ is the gradient
 operator and $\Delta$ is
Laplacian operator. The constant $\tau
 $ is the
 momentum relaxation time of electrons. The pressure
 $p$ satisfies the usual $\gamma$-law: $$p=p(n)=An^{\gamma},\eqno(1.2)$$
 where the case $\gamma>1$ corresponds to the isentropic gas and $\gamma=1$
 corresponds to the isothermal gas, $A$ is
 a positive
 constant.
 The constant $\bar{n}>0$ stands for the density of positively charged
 background
 ions. The system
 is supplemented with the initial data
 $$(n,\textbf{u})(x,0)=(n_{0},\textbf{u}_{0})(x),\ \
  x\in\mathbb{R}^{N}. \eqno(1.3)$$

 There are many contributions in mathematical analysis for (1.1)-(1.2), like
 well-posedness of steady state solutions, global existence of
 classical or entropy weak solutions, large time behavior of
 classical solutions, relaxation limit problems and so on, we may refer to
 Refs.~\cite{A,De,G,HMW,HW,LNX,WGHH} and the references therein.
For the evolutionary system (1.1)-(1.2), Luo, Natalini and Xin
\cite{LNX} first established
 the global exponential stability of small classical solutions to the Cauchy problem.
 Subsequently, there are many well-posedness and stability
 results
 for the Cauchy
 problem and the initial boundary value problem in one dimensional or
 multidimensional
 space. Physically, it is more important and more interesting to study (1.1)-(1.2) in the multidimensional case.
However, up to now, only partial results are available. Hsiao and
Wang \textit{et al.} \cite{HW,WGHH} studied the spherically
symmetrical solutions to (1.1)-(1.2) with $\gamma=1$, $\Omega=\{x\in
\mathbb{R}^{N}|\ 0<R_{1}\leqslant |x|\leqslant R_{2}<+\infty\}$ and
$\gamma>1$, $\Omega=\{x\in \mathbb{R}^{N}| \ |x|\geqslant R_{1}>0\}$
respectively. Guo \cite{G} investigated the irrotational
Euler-Poisson equation ($ N=3$) without relaxation and constructed
the global smooth irrotational solutions based on the Klein-Gordon
effect, which decay to the equilibrium state uniformly as
$(1+t)^{-p}(1<p<3/2)$. Recently, Hsiao, Markowich and Wang
\cite{HMW} dealt with
 the multidimensional unbounded domain problem ($N=2,3)$ without any
 geometrical assumptions. The main ingredient was to arrive at an \textit{a-priori} estimate
in terms of the \textit{classical energy argument}. Later,
Al\`{i}\cite{A} discussed the extended thermodynamic model
($N\geqslant2$) and reached the global existence, uniqueness and
stability of classical solutions. The key \textit{a-priori}
estimate was obtained by the positive definiteness of some
Liapunov functions. The above results are established on the
framework of Sobolev space $H^{\ell}(\mathbb{R}^{N})$ and the
regularity index is required to be high ($\ell>1+\frac{N}{2}, \
\ell\in \mathbb{Z}$) when one deals with them by classical
analysis methods. To lower the regularity, using the
Littlewood-Paley decomposition and Bony's para-product formula, we
prove the following well-posedness results for the system
(1.1)-(1.3) in the critical nonhomogeneous Besov space
$B^{\sigma}_{2,1}(\mathbb{R}^{N})$ ($\sigma=1+\frac{N}{2}$). First
of all, we give a local existence and uniqueness theorem of
classical solutions to (1.1)-(1.3) away from the vacuum.

\begin{thm} Let \ $\bar{n}>0$ be a constant reference density.
Suppose that\ $n_{0}-\bar{n}, \textbf{u}_{0}$ and
$\textbf{e}_{0}\in B^{\sigma}_{2,1}(\mathbb{R}^{N})$ with \
$n_{0}>0$, then there exist a time $T_{0}>0$ and a unique solution
$(n,\textbf{u},\nabla\Phi)$ of the system (1.1)-(1.3) such that
$$(n,\textbf{u},\nabla\Phi)\in \mathcal{C}^{1}([0,T_{0}]\times \mathbb{R}^{N})\ \ \ \mbox{with}\
\ \  n>0\ \  \mbox{for all} \ \ t\in [0,T_{0}]$$ and
$$(n-\bar{n},\textbf{u},\nabla\Phi)\in  \mathcal{C}([0,T_{0}],B^{\sigma}_{2,1}(\mathbb{R}^{N}))\cap
\mathcal{C}^1([0,T_{0}],B^{\sigma-1}_{2,1}(\mathbb{R}^{N})),$$
where $\textbf{e}_{0}:= \nabla\Delta^{-1}(n_{0}-\bar{n})$.
\end{thm}
\begin{rem}
 The symbol $\nabla\Delta^{-1}$ means
    $$\nabla\Delta^{-1}f=\int_{\mathbb{R}^N}\nabla_{x}G(x-y)f(y)dy,$$
 where $G(x,y)$ is a solution to $\Delta_xG(x,y)=\delta(x-y)$ with
 $x,y\in\mathbb{R}^N$.
\end{rem}
\begin{rem}
  (1) When one estimates  the nonlinear pressure
term by the spectral localization method, there appear many
commutators, which make computation very tedious. Thanks to the
ideas in Ref.~\cite{STW}, we introduce a function transform (sound
speed) to reduce the system (1.1)-(1.3) to a symmetric hyperbolic
system (3.1)-(3.2) where the nonlinear term becomes linear and
bi-linear terms. But, under the transform there is a new function
$h(m)$  in Poisson equation, which is well-defined and smooth on
the domain $\{m|\frac{\gamma-1}{2}m+\bar{\psi}>0\}$. Theorem 1.1
follows from Proposition 4.1 and Remark 3.1. We extend the local
well-posedness result in $H^{\ell}(\mathbb{R}^{N})$ (see
Ref.~\cite{A} or~\cite{HMW}) to that in Besov space
$B^{\sigma}_{2,1}(\mathbb{R}^{N})$ in virtue of the regularized
means and compactness argument.

 (2) Although our local existence result is proved via a symmetric
hyperbolic system, there are  some especial contents from Poisson
equation to be dealt with. Therefore, we can not use the result of
Iftimie \cite{I} directly.
\end{rem}

Under a smallness assumption, we establish the global existence,
uniqueness and exponential decay of classical solutions to
(1.1)-(1.3).
\begin{thm}  Let \ $\bar{n}>0$ be a constant reference density. Suppose that \ $n_{0}-\bar{n}, \textbf{u}_{0}$
and $\textbf{e}_{0}\in B^{\sigma}_{2,1}(\mathbb{R}^{N})$. There
exists a positive constant $\delta_{0}$ depending only on $A,
\gamma, \tau$ and $\bar{n}$ such that if
$$\|(n_{0}-\bar{n},\textbf{u}_{0},\textbf{e}_{0})\|_{B^{\sigma}_{2,1}(\mathbb{R}^{N})}
\leqslant \delta_{0},$$ then there exists a unique global solution
$(n,\textbf{u},\nabla\Phi)$ of the system (1.1)-(1.3) satisfying
$$(n,\textbf{u},\nabla\Phi)\in \mathcal{C}^{1}([0,\infty)\times \mathbb{R}^{N})$$
and
$$(n-\bar{n},\textbf{u},\nabla\Phi) \in
\mathcal{C}([0,\infty),B^{\sigma}_{2,1}(\mathbb{R}^{N}))\cap
\mathcal{C}^1([0,\infty),B^{\sigma-1}_{2,1}(\mathbb{R}^{N})).$$
Moreover, we have the decay estimate
\begin{eqnarray*}
&&\|(n-\bar{n},\textbf{u},\nabla\Phi)(\cdot,t)\|_{B^{\sigma}_{2,1}(\mathbb{R}^{N})}
+\|(n_{t},\textbf{u}_{t},\nabla\Phi_{t})(\cdot,t)\|_{B^{\sigma-1}_{2,1}(\mathbb{R}^{N})}
\nonumber\\&\leqslant&
C_{0}\|(n_{0}-\bar{n},\textbf{u}_{0},\textbf{e}_{0})\|_{B^{\sigma}_{2,1}(\mathbb{R}^{N})}\exp(-\mu_{0}t),\
 t\geqslant0,
\end{eqnarray*}
 where the positive
constants $\mu_{0}$ and $C_{0}$ depend only on $A, \gamma, \tau$
and $\bar{n}$.
\end{thm}
\begin{rem}
From the proof of Theorem 1.2, one can see that the decay rate
$\mu_{0}$ satisfies $\mu_{0}=K\tau$ if $0<\tau\leqslant1$ and
$\mu_{0}=K/\tau$ if $\tau>1$, where $K$ is a positive constant
independent of $\tau$.
\end{rem}

\begin{rem} Theorem 1.2 follows from Proposition 5.1 and Remark 3.1.
Proposition 5.1 relies on a crucial \textit{a-priori} estimate in
Proposition 5.2 and the standard boot-strap argument. In the proof
of \textit{a-priori} estimate, we use high- and low- frequency
decomposition method. Different from the \textit{classical energy
argument} in Ref.~\cite{HMW}, our method shows that Poisson
equation plays a key role in the low frequency estimates (see
(5.3), (5.4) and (5.10)). Such fact leads to the global
exponential stability of classical solutions to (3.1)-(3.2).
\end{rem}

Based on Theorem 1.2, we can characterize the exponential decay of
the vorticity $\omega=\nabla\times\textbf{u}$ in Besov space
$B^{\sigma-1}_{2,1}(\mathbb{R}^{N})$.

\begin{thm} (N=2,3) Let $(n,\textbf{u},\nabla\Phi)$ be the solution in Theorem 1.2. If
$$\|(n_{0}-\bar{n},\textbf{u}_{0},\textbf{e}_{0})\|_{B^{\sigma}_{2,1}(\mathbb{R}^{N})}
\leqslant \delta_{0}',$$ then the vorticity
$\omega=\nabla\times\textbf{u}$ decays exponentially in
$B^{\sigma-1}_{2,1}(\mathbb{R}^{N})(\omega(x,0)=\nabla\times
\textbf{u}_{0})$:
$$\|\omega(\cdot,t)\|_{B^{\sigma-1}_{2,1}(\mathbb{R}^{N})}\leqslant \|\omega(\cdot,0)\|_{B^{\sigma-1}_{2,1}(\mathbb{R}^{N})}\exp\Big(-\frac{t}{\tau}\Big), \ \ t\geqslant 0,$$
where the positive constant $\delta_{0}'=\min\{\delta_{0},
\frac{1}{2C_{5}C_{0}\tau}\}$ depends only on $A, \gamma, \tau$ and
$\bar{n}$ ($C_{5}$ a constant given in
 (5.22)).
\end{thm}

The paper is arranged as follows. In Section 2, we present some
definitions and basic facts on the Littlewood-Paley decomposition
and Bony's para-product formula. In Section 3, we reformulate the
system (1.1)-(1.3) in order to arrive at the effective
\textit{a-priori} estimates by the spectral localization method.
In Section 4, we obtain the local existence and uniqueness of
classical solutions to (3.1)-(3.2) with general initial data. In
the last section, we deduce a crucial \textit{a-priori} estimate
under a smallness assumption, which is used to achieve the proof
of global existence. Furthermore, it is also shown that the
vorticity decays to zero exponentially in the $2$D and $3$D space.

Throughout this paper, the symbol $C$ denotes a harmless constant
and  all functional spaces are considered in $\mathbb{R}^{N}$, so
we may omit the space dependence for simplicity. Moreover, the
integration $\int_{\mathbb{R}^{N}}fdx $ is labeled as $\int f$
without any ambiguity.

\section{Littlewood-Paley Analysis}
In this section, these definitions and basic facts can be found in
Darchin's \cite{D3} mini-course.

 Let $\mathcal{S}$ be the Schwarz class. ($\varphi, \chi)$ is a couple of smooth functions valued in [0,1]
such that $\varphi$ is supported in the shell
$\textbf{C}(0,\frac{3}{4},\frac{8}{3})=\{\xi\in\mathbb{R}^{N}|\frac{3}{4}\leqslant|\xi|\leqslant\frac{8}{3}\}$,
$\chi$ is supported in the ball $\textbf{B}(0,\frac{4}{3})=
\{\xi\in\mathbb{R}^{N}||\xi|\leqslant\frac{4}{3}\}$ and $$
\chi(\xi)+\sum_{q=0}^{\infty}\varphi(2^{-q}\xi)=1,\ \ \ \ q\in
\mathbb{Z},\ \  \xi\in\mathbb{R}^{N}.$$ For
$f\in\mathcal{S'}$(denote the set of temperate distributes which
is the dual one of $\mathcal{S}$), we can define the
nonhomogeneous dyadic blocks as follows:
$$\Delta_{-1}f:=\chi(D)f=\tilde{h}\ast f\ \ \ \mbox{with}\ \ \tilde{h}=\mathcal{F}^{-1}\chi,$$
$$\Delta_{q}f:=\varphi(2^{-q}D)f=2^{qN}\int h(2^{q}y)f(x-y)dy\ \ \ \mbox{with}\ \ h=\mathcal{F}^{-1}\varphi,\ \ \mbox{if}\ \ q\geqslant0.$$
where $\ast, \ \ \mathcal{F}^{-1} $ represent the convolution
operator and the inverse Fourier transform respectively. The
nonhomogeneous Littlewood-Paley decomposition is$$ f=\sum_{q
\geqslant-1}\Delta_{q}f \ \ \ \mbox{in}\ \ \ \mathcal{S'}.$$
Define the low frequency cut-off by
$$S_{q}f:=\sum_{p\leqslant q-1}\Delta_{p}f.$$ Of course, $S_{0}f=\Delta_{-1}f$. The above
Littlewood-Paley decomposition is almost orthogonal in $L^2$.
\begin{prop}
For any $f, g\in\mathcal{S'}$, the following properties hold:
$$\Delta_{p}\Delta_{q}f\equiv 0 \ \ \ \mbox{if}\ \ \ |p-q|\geqslant 2,$$
$$\Delta_{q}(S_{p-1}f\Delta_{p}g)\equiv 0\ \ \ \mbox{if}\ \ \ |p-q|\geqslant 5.$$
\end{prop}
 Besov space can be characterized in virtue of the
Littlewood-Paley decomposition.
\begin{defn}
Let $1\leqslant p\leqslant\infty$ and $s\in \mathbb{R}$. For
$1\leqslant r<\infty$, the Besov spaces  $B^{s}_{p,r}$ are defined
by
$$f\in B^{s}_{p,r} \Leftrightarrow \Big(\sum_{q\geqslant-1}(2^{qs}\|\Delta_{q}f\|_{L^{p}})^{r}\Big)^{\frac{1}{r}}<\infty$$
and $B^{s}_{p,\infty}$ are defined by
$$f\in B^{s}_{p,\infty} \Leftrightarrow \sup_{q\geqslant-1}2^{qs}\|\Delta_{q}f\|_{L^{p}}<\infty.$$
\end{defn}

\begin{defn}(J.-M. Bony \cite{B})
Let $f,g $ be two temperate distributions. The product $f\cdot g$
has the Bony's decomposition:
$$f\cdot g=\textbf{T}_{f}g+\textbf{T}_{g}f+\textbf{R}(f,g), $$
where $\textbf{T}_{f}g$ is paraproduct of $g$ by $f$,
$$ \textbf{T}_{f}g=\sum_{p\leqslant q-2}\Delta_{p}f\Delta_{q}g=\sum_{q}S_{q-1}f\Delta_{q}v$$
and the remainder $ \textbf{R}(f,g)$ is denoted by
$$\textbf{R}(f,g)=\sum_{q}\Delta_{q}f\tilde{\Delta}_{q}g\ \ \ \mbox{with} \ \
\tilde{\Delta}_{q}:=\Delta_{q-1}+\Delta_{q}+\Delta_{q+1}.$$
\end{defn}
\begin{lem}(Bernstein)
Let $k\in\mathbb{N}$ and $0<R_{1}<R_{2}$. There exists a constant
$C$ depending only on $R_{1},R_{2}$ and $N$ such that for all
$1\leqslant a\leqslant b\leqslant\infty$ and $f\in L^{a}$, we have
$$\mathrm{Supp}\ \mathcal{ F}f\subset \textbf{B}(0,R_{1}\lambda)\Rightarrow\sup_{|\alpha|=k}\|\partial^{\alpha}f\|_{L^{b}}\leqslant C^{k+1}\lambda^{k+N(\frac{1}{a}-\frac{1}{b})}\|f\|_{L^{a}};$$
$$\mathrm{Supp}\ \mathcal{ F}f\subset \textbf{C}(0,R_{1}\lambda,R_{2}\lambda)\Rightarrow C^{-k-1}\lambda^{k}\|f\|_{L^{a}}\leqslant \sup_{|\alpha|=k}\|\partial^{\alpha}f\|_{L^{a}}\leqslant C^{k+1}\lambda^{k}\|f\|_{L^{a}}.$$
Here, $\mathcal{F}$ represents the Fourier transform.
\end{lem}
A result of compactness in Besov space is:
\begin{prop}
Let $1\leqslant p,r\leqslant \infty,\ s\in \mathbb{R}$ and
$\varepsilon>0$. For all $\phi\in C_{c}^{\infty}$, the map
$f\mapsto\phi f$ is compact from $B^{s+\varepsilon}_{p,r}$ to
$B^{s}_{p,r}$.
\end{prop}
Finally, we state a result of continuity for the composition to
end this section.
\begin{prop}
Let $1\leqslant p,r\leqslant \infty$,\ \ $I$ be open interval of
$\mathbb{R}$. Let $s>0$ and $n$ be the smallest integer such that
$n\geqslant s$. Let $F:I\rightarrow\mathbb{R}$ satisfy $F(0)=0$
and $F'\in W^{n,\infty}(I;\mathbb{R}).$ Assume that $v\in
B^{s}_{p,r}$ takes values in $J\subset\subset I$. Then $F(v)\in
B^{s}_{p,r}$ and there exists a constant $C$ depending only on
$s,I,J$ and $N$ such that $$\|F(v)\|_{B^{s}_{p,r}}\leqslant
C(1+\|v\|_{L^{\infty}})^{n}\|F'\|_{W^{n,\infty}(I)}\|v\|_{B^{s}_{p,r}}.$$
\end{prop}
\section{Reformulation of the Original System}
In this section, we are going to reformulate (1.1)-(1.3) in order
to obtain the effective \textit{a-priori} estimates by spectral
localization method. For the isentropic case $(\gamma>1)$, let the
sound speed
$$\psi(n)=\sqrt{p'(n)},$$ and denote the sound speed at a background density $\bar{n}$ by  $\bar{\psi}=\psi(\bar{n})$. Similar to that in Ref.~\cite{STW}, we define
$$m=\frac{2}{\gamma-1}\Big(\psi(n)-\bar{\psi}\Big).$$
Then the system (1.1) can be reduced to the following system for
$\mathcal{C}^1$ solutions:
$$\begin{cases}
m_{t}+\bar{\psi}\mbox{div}\textbf{u}=-\textbf{u}\cdot\nabla
m-\frac{\gamma-1}{2}m\mbox{div}\textbf{u},\cr
 \textbf{u}_{t}+\bar{\psi}\nabla m+\frac{\textbf{u}}{\tau}=-\textbf{u}\cdot\nabla\textbf{u}-\frac{\gamma-1}{2}m\nabla m+\textbf{e}, \cr
 \textbf{e}_{t}=-\nabla\Delta^{-1}\nabla\cdot\{h(m)\textbf{u}+\bar{n}\textbf{u}\},
 \end{cases}\eqno(3.1)$$
where $\textbf{e}:=\nabla\Phi$,
$h(m)=\{(A\gamma)^{-\frac{1}{2}}(\frac{\gamma-1}{2}m+\bar{\psi})\}^{\frac{2}{\gamma-1}}-\bar{n}$\
is a smooth function on the domain
$\{m|\frac{\gamma-1}{2}m+\bar{\psi}>0\}$ satisfying $h(0)=0$, the
non-local term $\nabla\Delta^{-1}\nabla \cdot f$ is the product of
Riesz transforms of $f$. The initial data (1.3) becomes into
$$(m,\textbf{u},\textbf{e})|_{t=0}=(m_{0},\textbf{u}_{0},\textbf{e}_{0})\eqno(3.2)$$
with $$m_{0}=\frac{2}{\gamma-1}\Big(\psi(n_{0})-\bar{\psi}\Big), \
\ \ \textbf{e}_{0}=\nabla\Delta^{-1}(n_{0}-\bar{n}).$$

\begin{rem} For any $T>0$,
$(n,\textbf{u},\textbf{e})\in \mathcal{C}^1([0,T]\times
\mathbb{R}^{N})$ is a solution of the system (1.1)-(1.2) with
$n>0$, then $(m,\textbf{u},\textbf{e})\in
\mathcal{C}^1([0,T]\times \mathbb{R}^{N})$ solves the system
(3.1)-(3.2) with $\frac{\gamma-1}{2}m+\bar{\psi}>0$; Conversely,
if $(m,\textbf{u},\textbf{e})\in \mathcal{C}^1([0,T]\times
\mathbb{R}^{N})$ solves the system (3.1)-(3.2) with
$\frac{\gamma-1}{2}m+\bar{\psi}>0$, then
$(n,\textbf{u},\textbf{e})\in \mathcal{C}^1([0,T]\times
\mathbb{R}^{N})$ is a solution of the system (1.1)-(1.2) with
$n>0$, where $n=\psi^{-1}(\frac{\gamma-1}{2}m+\bar{\psi})$.
\end{rem}
 For the isothermal case $(\gamma=1)$, let $\tilde{n}=\sqrt{A}(\ln n-\ln
\bar{n}),\  \textbf{e}=\nabla\Phi$, then the system (1.1) is
transformed into the following one for $\mathcal{C}^1$ solutions:
$$\begin{cases}\tilde{n}_{t}+\sqrt{A}\mbox{div}\textbf{u}=-\textbf{u}\cdot\nabla \tilde{n},\cr
 \textbf{u}_{t}+\sqrt{A}\nabla \tilde{n}+\frac{\textbf{u}}{\tau}=-\textbf{u}\cdot\nabla\textbf{u}+\textbf{e}, \cr
 \textbf{e}_{t}
 =-\nabla\Delta^{-1}\nabla\cdot\{h(\tilde{n})\textbf{u}+\bar{n}\textbf{u}\},
 \end{cases}\eqno(3.3)$$
where $h(\tilde{n})=\bar{n}(\exp(A^{-\frac{1}{2}}\tilde{n})-1)$ is
a smooth  function on the domain
$\{\tilde{n}|-\infty<\tilde{n}<+\infty\}$ satisfying $h(0)=0$. The
initial data (1.3) turns into
$$(\tilde{n},\textbf{u},\textbf{e})|_{t=0}=(\sqrt{A}(\ln n_{0}-\ln
\bar{n}), \textbf{u}_{0},\textbf{e}_{0}). \eqno(3.4) $$
\begin{rem} For any $T>0$, if $(n,\textbf{u},\textbf{e})\in
\mathcal{C}^1([0,T]\times \mathbb{R}^{N})$ \ is a solution of the
system (1.1)-(1.2) with $n>0$, \ then
$(\tilde{n},\textbf{u},\textbf{e})\in \mathcal{C}^1([0,T]\times
\mathbb{R}^{N})$ solves the system (3.3)-(3.4); Conversely, if
$(\tilde{n},\textbf{u},\textbf{e})\in \mathcal{C}^1([0,T]\times
\mathbb{R}^{N})$ solves the system (3.3)-(3.4), then
$(n,\textbf{u},\textbf{e})\in \mathcal{C}^1([0,T]\times
\mathbb{R}^{N})$ is a solution of the system (1.1)-(1.2) with
$n>0$, where $n=\bar{n}\exp(A^{-\frac{1}{2}}\tilde{n})$.
\end{rem}
In what follows, we shall only study the system (3.1)-(3.2) and
prove the main results in this paper, since (3.3)-(3.4) can be
discussed through a similar process.
\section{Local Existence}
In this section, we first give some  estimates of commutators in
Besov space $B^{s}_{p,1}$. Then, using the regularized means and
compactness argument, we complete the proof of proposition 4.1.

By using the first and third equation of Eq.(1.1), we get the
following formulations under the variable transform immediately.
\begin{lem}
$$\mathrm{div}\textbf{e}=h(m),\ \ \ \mathrm{div}\textbf{u}=-\frac{\mathrm{div}\textbf{e}_{t}+\mathrm{div}(h(m)\textbf{u})}{\bar{n}}.$$
\end{lem}
\hspace{118mm}$\square$

Applying the operator $\Delta_{q}$ to (3.1) yields
\begin{equation}
\left\{
\begin{array}{l}
\partial_{t}\Delta_{q}m+(\textbf{u}\cdot\nabla)\Delta_{q}m=-\bar{\psi}\Delta_{q}\mbox{div}\textbf{u}+[\textbf{u},\Delta_{q}]\cdot\nabla
m-\frac{\gamma-1}{2}\Delta_{q}(m\mbox{div}\textbf{u}),\\
 \partial_{t}\Delta_{q}\textbf{u}+(\textbf{u}\cdot\nabla)\Delta_{q}\textbf{u}+\frac{1}{\tau}\Delta_{q}\textbf{u}
 =-\bar{\psi}\Delta_{q}(\nabla m)+[\textbf{u},\Delta_{q}]\cdot\nabla\textbf{u}
 \\\hspace{53mm}-\frac{\gamma-1}{2}\Delta_{q}(m\nabla
 m)+\Delta_{q}\textbf{e},\\
 \partial_{t}\Delta_{q}\textbf{e}=-\nabla\Delta^{-1}\nabla\cdot\Delta_{q}\{h(m)\textbf{u}+\bar{n}\textbf{u}\},
\end{array}
\right.
\end{equation}
where the commutator $[f,g]=fg-gf.$ \\  Multiplying the first
equation of Eq.(4.1) by $\Delta_{q}m$, the second one by
$\Delta_{q}\textbf{u}$ and adding the resulting equations
together, then integrating it over $\mathbb{R}^{N}$, we obtain
\begin{eqnarray}&&\frac{1}{2}\frac{d}{dt}\Big(\|\Delta_{q}m\|^2_{L^2}+\|\Delta_{q}\textbf{u}\|^2_{L^2}\Big)+\frac{1}{\tau}\|\Delta_{q}\textbf{u}\|^2_{L^2}\nonumber
\\&=&\frac{1}{2}\int\mathrm{div}\textbf{u}(|\Delta_{q}m|^2+|\Delta_{q}\textbf{u}|^2)+\int([\textbf{u},\Delta_{q}]\cdot\nabla
m\Delta_{q}m+[\textbf{u},\Delta_{q}]\cdot\nabla\textbf{u}\Delta_{q}\textbf{u})
\nonumber\\&&-\frac{\gamma-1}{2}\int\Big(\Delta_{q}(m\mathrm{div}\textbf{u})\Delta_{q}m+\Delta_{q}(m\nabla
m)\Delta_{q}\textbf{u}\Big)+\int\Delta_{q}\textbf{e}\Delta_{q}\textbf{u}.\end{eqnarray}
The electric field term can be estimated as
\begin{eqnarray}&&\int\Delta_{q}\textbf{e}\Delta_{q}\textbf{u}\nonumber\\&=&-\int\Delta_{q}\Phi\Delta_{q}\mathrm{div}\textbf{u}\nonumber\\
&=&\frac{1}{\bar{n}}\int\Delta_{q}\Phi\Delta_{q}\Big(\mathrm{div}\textbf{e}_{t}+\mathrm{div}(h(m)\textbf{u})\Big)\nonumber\\&=&-\frac{1}{2\bar{n}}\frac{d}{dt}\|\Delta_{q}\textbf{e}\|_{L^2}^2-\frac{1}{\bar{n}}\int\Delta_{q}\textbf{e}\Delta_{q}(h(m)\textbf{u}).\end{eqnarray}
Note that the bi-linear spectral localization term, we have
\begin{eqnarray}
&&\int\{\Delta_{q}(m\mathrm{div}\textbf{u})\Delta_{q}m+\Delta_{q}(m\nabla
m)\Delta_{q}\textbf{u}\}\nonumber\\&=&-\int\Delta_{q}m(\nabla
m\cdot\Delta_{q}\textbf{u})+\int[\Delta_{q},m]\nabla
m\cdot\Delta_{q}\textbf{u}+\int[\Delta_{q},m]\mathrm{div}\textbf{u}\Delta_{q}m.\end{eqnarray}
Here, we give a lemma to estimate these commutators in (4.2) and
(4.4).
\begin{lem} Let $s>0$ and $1<p<\infty$, the following inequalities are true:
\begin{eqnarray}
&&2^{qs}\|[f,\Delta_{q}]\mathcal{A}
g\|_{L^{p}}\nonumber\\&\leqslant& \left\{
\begin{array}{l}
 Cc_{q}\|f\|_{B^{s}_{p,1}}\|g\|_{B^{s}_{p,1}},\ \ \
\ \ \ f,g\in B^{s}_{p,1},\ s=1+N/p\ ;\\
 Cc_{q}\|f\|_{B^{s}_{p,1}}\|g\|_{B^{s+1}_{p,1}},\
 \ \ \ f\in B^{s}_{p,1},\ \ g\in B^{s+1}_{p,1}, \ s=N/p\ ;\\
 Cc_{q}\|f\|_{B^{s+1}_{p,1}}\|g\|_{B^{s}_{p,1}},\ \
\ \ \ f\in B^{s+1}_{p,1},\ \ g\in B^{s}_{p,1}, \ s=N/p\ .
\end{array} \right.
\end{eqnarray}
In particular, if $f=g$, then
\begin{equation}
2^{qs}\|[f,\Delta_{q}]\mathcal{A}g\|_{L^{p}} \leqslant
Cc_{q}\|\nabla f\|_{L^{\infty}}\|g\|_{B^{s}_{p,1}},\ s>0.
\end{equation}
where the operator $\mathcal{A}=\mathrm{div}$ or
$\mathrm{\nabla}$, $C$ is a harmless  constant and $\{c_{q}\}$
denotes a sequence such that $\|\{c_{q}\}\|_{ {l^{1}}}\leqslant
1.$
\end{lem}
\begin{rem}
 Similar estimates of commutators  have been obtained by
  Chemin\cite{C}, Danchin\cite{D1,D2,D3} and Iftimie\cite{I} \textit{et
  al.},
 therefore, we omit the details.
\end{rem}
Now, we give a proposition on the local existence and uniqueness
of classical solutions to (3.1)-(3.2).
\begin{prop}
Suppose that
$(m_{0},\textbf{u}_{0},\textbf{e}_{0})\in{B^{\sigma}_{2,1}}$
satisfying $ \frac{\gamma-1}{2}m_{0}+\bar{\psi}>0$, then there
exist a time $T_{0}>0$ and a unique solution $(m, \textbf{u},
\textbf{e})$ of (3.1)-(3.2) such that $(m, \textbf{u},
\textbf{e})\in \mathcal{C}^{1}([0,T_{0}]\times \mathbb{R}^{N})$
with $ \frac{\gamma-1}{2}m+\bar{\psi}>0$ for \mbox{all}
$t\in[0,T_{0}]$ and $(m, \textbf{u}, \textbf{e})\in
\mathcal{C}([0,T_{0}],B^{\sigma}_{2,1})\cap
\mathcal{C}^1([0,T_{0}],B^{\sigma-1}_{2,1})$.
\end{prop}
\begin{proof} (\textbf{Existence}) Let $U_{0}=(m_{0},
\textbf{u}_{0},\textbf{e}_{0})^{\top}\in{B^{\sigma}_{2,1}}$.

\noindent\textbf{Claim 1}: There exist two positive constants
$\eta', \eta''$ such that
\begin{eqnarray} 0<\eta'\leqslant\frac{\gamma-1}{2}m_{0}+\bar{\psi}\leqslant\eta''.\end{eqnarray}
In fact, by the imbedding $B^{\sigma-1}_{2,1}\hookrightarrow
\mathcal{C}_{0}$(continuous bounded functions which decay to zero
at infinity), we know $m_{0}\in \mathcal{C}_{0}$.  For any
$\varepsilon>0$, $\exists M=M(\varepsilon)>0$ such that
$|m_{0}|\leqslant\varepsilon$ for all $|x|>M$. We may choose
$\varepsilon=\frac{1}{\gamma-1}\bar{\psi}$, then
$0<\frac{\bar{\psi}}{2}\leqslant\frac{\gamma-1}{2}m_{0}+\bar{\psi}\leqslant\frac{3\bar{\psi}}{2}$.
On the other hand, when $x\in[-M,M]$, we have
$$0<\frac{\gamma-1}{2}\min_{x\in[-M,M]} m_{0}+\bar{\psi}\leqslant\frac{\gamma-1}{2}m_{0}+\bar{\psi}\leqslant\frac{\gamma-1}{2}\max_{x\in[-M,M]} m_{0}+\bar{\psi}.$$
So we obtain (4.7) only by choosing
$$\eta'=\max\Big\{\frac{\bar{\psi}}{2},
\frac{\gamma-1}{2}\min_{x\in[-M,M]} m_{0}+\bar{\psi}\Big\},\ \ \ \
\ \eta''=\min\Big\{\frac{3\bar{\psi}}{2},
\frac{\gamma-1}{2}\max_{x\in[-M,M]} m_{0}+\bar{\psi}\Big\} .$$

From (4.7), we can obtain $0<\vartheta'\leqslant n_{0}\leqslant
\vartheta''$ for two positive constants $\vartheta', \vartheta''$.
Note that
$$h(m_{0})=\Big\{(A\gamma)^{-\frac{1}{2}}\Big(\frac{\gamma-1}{2}m_{0}+\bar{\psi}\Big)\Big\}^{\frac{2}{\gamma-1}}-\bar{n}=n_{0}-\bar{n}$$ is a
smooth function on $
\{m_{0}|\frac{\gamma-1}{2}m_{0}+\bar{\psi}>0\}$, by Proposition
2.3, we have
$$\|n_{0}-\bar{n}\|_{B^{\sigma}_{2,1}}=\|h(m_{0})\|_{B^{\sigma}_{2,1}}\leqslant C\|m_{0}\|_{B^{\sigma}_{2,1}}.$$
There exists a sequence $\{\tilde{U}^{k}_{0}\}:=\{(n_{0}^{k},
\textbf{u}^{k}_{0}, \textbf{e}^{k}_{0})^{\top}\}$ such that
$\{(n_{0}^{k}-\bar{n}, \textbf{u}^{k}_{0},
\textbf{e}^{k}_{0})\}\in H^{\ell} (\ell>\sigma, \
\ell\in\mathbb{Z})$ converges to $(n_{0}-\bar{n}, \textbf{u}_{0},
\textbf{e}_{0})$ in $B^{\sigma}_{2,1}$\ and \ satisfies
$0<\frac{1}{2}\vartheta'\leqslant n_{0}^{k}\leqslant
\vartheta''+\frac{1}{2}\vartheta'$.\ \ Furthermore, we also have
$\|U^{k}_{0}\|_{B^{\sigma}_{2,1}}\leqslant
\|U_{0}\|_{B^{\sigma}_{2,1}}+1 (U^{k}_{0}:=(m^{k}_{0},
\textbf{u}^{k}_{0},\textbf{e}^{k}_{0})^{\top})$. We define a
sequence $\{\tilde{U}^{k}\}=\{(n^{k}, \textbf{u}^{k},
\textbf{e}^{k})^{\top}\}$, which solves the following equations:
\begin{equation}\left\{
\begin{array}{l}n^{k}_{t}+\nabla\cdot(n^{k}\textbf{u}^{k})=0\\
 \textbf{u}^{k}_{t}+(\textbf{u}^{k}\cdot\nabla)\textbf{u}^{k}+\frac{1}{n^{k}}\nabla p(n^{k})=\textbf{e}^{k}-\frac{\textbf{u}^{k}}{\tau}
 \\
 \textbf{e}_{t}^{k}=-\nabla\Delta^{-1}\nabla\cdot(n^{k}\textbf{u}^{k})\end{array}
\right.\end{equation} with the initial data
\begin{equation}(
n^{k},\textbf{u}^{k},\textbf{e}^{k})|_{t=0}=(n^{k}_{0},\textbf{u}^{k}_{0},\textbf{e}^{k}_{0}).\end{equation}
It is easy to see (4.8) is a strict hyperbolic symmetric system on
$\tilde{G}=\{\tilde{U}^{k}|n^{k}>0\}$ in the sense of Friedrichs.
Using Kato's classical result in Ref.~\cite{K} or \cite{M} and the
$L^2$- boundedness of Riesz transformation, we can get the
following local existence result: there exist a time $T_{k}>0$ and
a solution $\tilde{U}^{k}$ to (4.8)-(4.9) such
that$$\tilde{U}^{k}\in \mathcal{C}^{1}([0,T_{k}]\times
\mathbb{R}^{N})\ \ \ \mbox{with}\ \ \ 0<\vartheta'_{k}\leqslant
n^{k}\leqslant\vartheta''_{k}\ \ \mbox{for all}\ t\in[0,T_{k}]$$$$
(\vartheta'_{k},\vartheta''_{k}\ \mbox{are two positive constants
depending on $k$})$$ and
$$(n^{k}-\bar{n},\textbf{u}^{k}, \textbf{e}^{k})\ \in \mathcal{C}([0,T_{k}],H^{\ell})\cap
\mathcal{C}^{1}([0,T_{k}],H^{\ell-1}).$$
 From Remark
3.1, the sequence $\{U^{k}\}=\{(m^{k}, \textbf{u}^{k},
\textbf{e}^{k})^{\top}\}$ solves the following equations:
\begin{equation}\left\{
\begin{array}{l} m^{k}_{t}+\bar{\psi}\mbox{div}\textbf{u}^{k}=-\textbf{u}^{k}\cdot\nabla
m^{k}-\frac{\gamma-1}{2}m^{k}\mbox{div}\textbf{u}^{k}\\
 \textbf{u}^{k}_{t}+\bar{\psi}\nabla m^{k}+\frac{\textbf{u}^{k}}{\tau}=-\textbf{u}^{k}\cdot\nabla\textbf{u}^{k}-\frac{\gamma-1}{2}m^{k}\nabla m^{k}+\textbf{e}^{k}
 \\
 \textbf{e}^{k}_{t}=-\nabla\Delta^{-1}\nabla\cdot\{h(m^{k})\textbf{u}^{k}+\bar{n}\textbf{u}^{k}\}
\end{array}
\right.\end{equation} with the initial data
\begin{equation}(
m^{k},\textbf{u}^{k},\textbf{e}^{k})|_{t=0}=\Big(\frac{2}{\gamma-1}\Big(\psi(n^{k}_{0})-\bar{\psi}\Big),\textbf{u}^{k}_{0},\textbf{e}^{k}_{0}\Big)\end{equation}
satisfying
$$U^{k}\in \mathcal{C}^{1}([0,T_{k}]\times \mathbb{R}^{N})\ \ \ \mbox{with}\ \  0<\eta'_{k}\leqslant\frac{\gamma-1}{2}m^{k}+\bar{\psi}\leqslant \eta''_{k}\ \ \mbox{for all}\ t\in[0,T_{k}]$$$$ (\eta'_{k},\eta''_{k} \ \mbox{are two positive
constants depending on $k$})$$ and $$U^{k}\in
\mathcal{C}([0,T_{k}],H^{\ell})\cap
\mathcal{C}^{1}([0,T_{k}],H^{\ell-1}).$$ Let $[0,T^{*}_{k})$ be
the maximal existence interval of above solutions to
(4.10)-(4.11). Similar to the discussion in Ref.~\cite{M}, we have
the blow-up criterion:
\begin{eqnarray*}T^{*}_{k}<\infty &\Leftrightarrow& \limsup_{t\rightarrow T^{*}_{k}}(\|U_{t}^{k}\|_{L^{\infty}}+\|\nabla
{U}^{k}\|_{L^{\infty}})=+\infty\mbox{\ \ \ or}\\
&&\mbox{for any compact subset} \ \ K \subset\subset G, \ \
U^{k}(x,t)\ \ \mbox{escapes}\ \  K \ \ \mbox{as} \ \ t\rightarrow
T^{*}_{k}, \\ && \mbox{where}\
G=\Big\{U^{k}|\frac{\gamma-1}{2}m^{k}+\bar{\psi}>0\Big\}.
\end{eqnarray*}
\noindent\textbf{Claim 2}: For $t\in [0,\min\{T^{*}_{k},T_{0}\}),$
it yields $\|U^{k}(t)\|_{B^{\sigma}_{2,1}}\leqslant 4\lambda_{0},$
where $T_{0}=\frac{1}{2\tilde{C}\lambda_{0}}$,
$\lambda_{0}=\|U_{0}\|_{B^{\sigma}_{2,1}}+1$
 and $\tilde{C}$ is a positive
constant (independent of $k$) given in (4.18).

The proof of Claim 2 relies on the standard continuity method. Let
\begin{equation}\mathcal{G}=\Big\{T\in[0,\min\{T^{*}_{k},T_{0}\}):\sup_{0\leqslant t\leqslant T}\mathcal{E}(t)\leqslant4\lambda_{0}\Big\},\end{equation}
where$$\mathcal{E}(t)=\|U^{k}(t)\|_{B^{\sigma}_{2,1}}.$$ In fact,
we have already known
\begin{eqnarray*}\|U^{k}_{0}\|_{B^{\sigma}_{2,1}}\leqslant
\lambda_{0}\leqslant 4\lambda_{0}.\end{eqnarray*} Then from the
continuity of $\mathcal{E}(t)$ on $[0,\min\{T^{*}_{k},T_{0}\})$,
we may see that the set $\mathcal{G}$ is nonempty and relatively
close in $[0,\min\{T^{*}_{k},T_{0}\})$. To show that it is also
relatively open in $[0,\min\{T^{*}_{k},T_{0}\})$, and hence the
entire interval $[0,\min\{T^{*}_{k},T_{0}\})$, it suffices to the
weaker bound in (4.12) implies
\begin{eqnarray}
\|U^{k}(t)\|_{B^{\sigma}_{2,1}}\leqslant2\lambda_{0},\ \ \ t\in
[0,T]\subset[0,\min\{T^{*}_{k},T_{0}\}).\end{eqnarray}

 To do this, we need to make the best use of
flow map.

 The flow map
$X(t;\bar{t},x)$ of $\textbf{u}^{k}$ starting from $x\in
\mathbb{R}^{N}$ at time $\bar{t} \in[0,T]$ can be defined as
$$\frac{dX}{dt}=\textbf{u}^{k}(t,X(t;\bar{t},x)),\ \ \ X(t;\bar{t},x)|_{t=\bar{t}}=x,$$
then we have
$$\frac{d}{dt}m^{k}(t,X(t;\bar{t},x))=(\partial_{t}+\textbf{u}^{k}\cdot\nabla)m^{k}(t,X(t;\bar{t},x)).$$
Together with the first equation of Eq.(4.10), it is easy
 to get (for any $\bar{t} \in[0,T]$)
\begin{eqnarray}
\hspace{-5mm}&&\frac{\gamma-1}{2}m^{k}(\bar{t},x)+\bar{\psi}\nonumber\\
\hspace{-5mm}&=&\Big(\frac{\gamma-1}{2}m^{k}(0,X(0;\bar{t},x))+\bar{\psi}\Big)\exp\Big(
-\frac{\gamma-1}{2}\int^{\bar{t}}_{0}\mathrm{div}\textbf{u}^{k}(\varsigma,X(\varsigma;\bar{t},x))d\varsigma\Big).
\end{eqnarray}
Thus, there exist two positive constants $\eta_{1},\eta_{2}\
$(independent of $k$) such that
\begin{eqnarray}0<\eta_{1}\leqslant\frac{\gamma-1}{2}m^{k}(\bar{t},x)+\bar{\psi}\leqslant\eta_{2},\
\ \ \ (\bar{t},x)\in[0,T]\times\mathbb{R}^{N}.\end{eqnarray} Here,
we need not consider the effect of relaxation term. Therefore, by
(4.2)-(4.4), H\"{o}lder's inequality and Lemma 4.2 ( take $p=2$
and $s=\sigma$), we get
\begin{eqnarray}&&\frac{d}{dt}\Big(2^{2q\sigma}\|\Delta_{q}U^{k}\|^2_{L^2}\Big)\nonumber\\\hspace{-5mm}&\leqslant&
C\{\|\nabla\textbf{u}^{k}\|_{L^{\infty}}2^{2q\sigma}(\|\Delta_{q}m^{k}\|^2_{L^2}+\|\Delta_{q}\textbf{u}^{k}\|^2_{L^2})+c_{q}2^{q\sigma}\|\textbf{u}^{k}\|_{B^{\sigma}_{2,1}}\|m^{k}\|_{B^{\sigma}_{2,1}}\|\Delta_{q}m^{k}\|_{L^2}\nonumber\\
&&+c_{q}2^{q\sigma}\|\nabla
\textbf{u}^{k}\|_{L^{\infty}}\|\textbf{u}^{k}\|_{B^{\sigma}_{2,1}}\|\Delta_{q}\textbf{u}^{k}\|_{L^2}+2^{2q\sigma}\|\nabla
m^{k}\|_{L^{\infty}}\|\Delta_{q}\textbf{u}^{k}\|_{L^2}\|\Delta_{q}m^{k}\|_{L^2}\nonumber\\&&+c_{q}2^{q\sigma}\|\nabla
m^{k}\|_{L^{\infty}}\|m^{k}\|_{B^{\sigma}_{2,1}}\|\Delta_{q}\textbf{u}^{k}\|_{L^2}+c_{q}2^{q\sigma}\|m^{k}\|_{B^{\sigma}_{2,1}}\|\textbf{u}^{k}\|_{B^{\sigma}_{2,1}}\|\Delta_{q}m^{k}\|_{L^2}
\nonumber\\&&+2^{2q\sigma}\|\Delta_{q}(h(m^{k})\textbf{u}^{k})\|_{L^2}\|\Delta_{q}\textbf{e}^{k}\|_{L^2}\},\
\ \ \ t\in[0,T].
\end{eqnarray}
$$(\|\Delta_{q}U^{k}\|^2_{L^2}:=\|\Delta_{q}m^{k}\|^2_{L^2}+\|\Delta_{q}\textbf{u}^{k}\|^2_{L^2}+\|\Delta_{q}\textbf{e}^{k}\|^2_{L^2},$$$$\|\Delta_{q}U^{k}\|_{L^2}\thickapprox\|\Delta_{q}m^{k}\|_{L^2}+\|\Delta_{q}\textbf{u}^{k}\|_{L^2}+\|\Delta_{q}\textbf{e}^{k}\|_{L^2})$$\\
Dividing (4.16) by
$(2^{2q\sigma}\|\Delta_{q}U^{k}\|^2_{L^2}+\epsilon)^{1/2}$
($\epsilon>0$ a small quantity), we have
\begin{eqnarray}\hspace{-5mm}&&\frac{d}{dt}\Big(2^{2q\sigma}\|\Delta_{q}U^{k}\|^2_{L^2}+\epsilon\Big)^{\frac{1}{2}}\nonumber\\\hspace{-5mm}&\leqslant&
C\{\|\nabla\textbf{u}^{k}\|_{L^{\infty}}2^{q\sigma}(\|\Delta_{q}m^{k}\|_{L^2}+\|\Delta_{q}\textbf{u}^{k}\|_{L^2})+c_{q}\|\textbf{u}^{k}\|_{B^{\sigma}_{2,1}}\|m^{k}\|_{B^{\sigma}_{2,1}}
\nonumber\\\hspace{-5mm} &&+c_{q}\|\nabla
\textbf{u}^{k}\|_{L^{\infty}}\|\textbf{u}^{k}\|_{B^{\sigma}_{2,1}}+\|\nabla
m^{k}\|_{L^{\infty}}2^{q\sigma}\|\Delta_{q}\textbf{u}^{k}\|_{L^2}+c_{q}\|\nabla
m^{k}\|_{L^{\infty}}\nonumber\\\hspace{-5mm}&&\times\|m^{k}\|_{B^{\sigma}_{2,1}}+c_{q}\|m^{k}\|_{B^{\sigma}_{2,1}}\|\textbf{u}^{k}\|_{B^{\sigma}_{2,1}}+2^{q\sigma}\|\Delta_{q}(h(m^{k})\textbf{u}^{k})\|_{L^2}\},\
 t\in[0,T].
\end{eqnarray} Integrating (4.17) on the variable $t$, then taking
$\epsilon\rightarrow 0$ and using Proposition 2.3, we obtain the
\textit{a-priori} estimate of $U^{k}$:
\begin{eqnarray}\|U^{k}(t)\|_{B^{\sigma}_{2,1}}&\leqslant&\|U^{k}_{0}\|_{B^{\sigma}_{2,1}}+\tilde{C}\int^{t}_{0}\|U^{k}(\varsigma)\|_{B^{\sigma}_{2,1}}^2d\varsigma\nonumber\\
&\leqslant&(\|U_{0}\|_{B^{\sigma}_{2,1}}+1)+\tilde{C}\int^{t}_{0}\|U^{k}(\varsigma)\|_{B^{\sigma}_{2,1}}^2d\varsigma,
\ \ \ t\in[0,T].\end{eqnarray} Furthermore, we have
\begin{equation}\sup_{0\leqslant\varsigma\leqslant
t}\|U^{k}(\varsigma)\|_{B^{\sigma}_{2,1}}\leqslant
(\|U_{0}\|_{B^{\sigma}_{2,1}}+1)+\tilde{C}\int^{t}_{0}\sup_{0\leqslant\varsigma'\leqslant
\varsigma}\|U^{k}(\varsigma')\|_{B^{\sigma}_{2,1}}^2d\varsigma,\ \
\ t\in[0,T]. \end{equation} Set
$$\lambda_{1}(t)\equiv (\|U_{0}\|_{B^{\sigma}_{2,1}}+1)+\tilde{C}\int^{t}_{0}\sup_{0\leqslant\varsigma'\leqslant
\varsigma}\|U_{k}(\varsigma')\|_{B^{\sigma}_{2,1}}^2d\varsigma.$$
Then, we have
\begin{eqnarray}\frac{d}{dt}\lambda_{1}\leqslant \tilde{C}\lambda_{1}^2,\ \ \ \ \
\lambda_{1}(0)=\|U_{0}\|_{B^{\sigma}_{2,1}}+1,\ \
t\in[0,T].\end{eqnarray} Let $\lambda(t)$ solves Riccati equation:
\begin{eqnarray}\frac{d}{dt}\lambda=\tilde{C}\lambda^2,\ \ \ \ \
\lambda(0)=\|U_{0}\|_{B^{\sigma}_{2,1}}+1.\end{eqnarray} The time
$T_{0}=1/(2\tilde{C}\lambda_{0})$ is less than the blow-up time
for (4.21). Then by solving the differential inequality (4.20), we
have $\lambda_{1}(t)\leqslant\lambda(t)$ for $t\in [0,T_{0}].$
Solving (4.21) yields
\begin{eqnarray}\lambda_{1}(t)\leqslant\frac{\|U_{0}\|_{B^{\sigma}_{2,1}}+1}{1-\tilde{C}t(\|U_{0}\|_{B^{\sigma}_{2,1}}+1)}=\lambda(t),
\ \ \ t\in[0,T_{0}]\ .\end{eqnarray} Therefore, we get
\begin{eqnarray*}
\|U^{k}(t)\|_{B^{\sigma}_{2,1}}\leqslant2\lambda_{0}, \ \ \
t\in[0,T]\subset[0,\min\{T^{*}_{k},T_{0}\}),\end{eqnarray*}
 which completes the proof of Claim 2.

 Furthermore, using Eq.(4.10), we have
\begin{eqnarray}\|U^{k}_{t}(t)\|_{B^{\sigma-1}_{2,1}}\leqslant \lambda_{0}',\ \ \ t\in[0,\min\{T^{*}_{k},T_{0}\}),\end{eqnarray}
where $\lambda_{0}'$ is a positive constant only depending on the
initial data $U_{0}$. There exist two positive constants
$\hat{\eta}, \check{\eta} $ (independent of $k$) such that
\begin{eqnarray}0<\hat{\eta}\leqslant\frac{\gamma-1}{2}m^{k}(t,x)+\bar{\psi}\leqslant \check{\eta},\
\ \ \
(t,x)\in[0,\min\{T^{*}_{k},T_{0}\})\times\mathbb{R}^{N}.\end{eqnarray}
From Claim 2, (4.23) and (4.24), the blow-up criterion implies
$0<T_{0}<T^{*}_{k}$, moreover, we have $0<T_{0}\leqslant
\inf_{k}T^{*}_{k}$.

That is, we find a positive time $T_{0}$ (only depending on the
initial data $U_{0}$) such that the approximative solution
sequence $\{U^{k}\}$ to (4.10)-(4.11) is uniformly bounded in
$\mathcal{C}([0,T_{0}],B^{\sigma}_{2,1})\cap
\mathcal{C}^1([0,T_{0}],B^{\sigma-1}_{2,1}).$ Moreover, it
weak$^{\star}$-converges (up to a subsequence) to some $U$ in
$L^{\infty}([0,T_{0}],B^{\sigma}_{2,1})$ in terms of the
Banach-Alaoglu Theorem (see Ref.~\cite{T} Remark 2 on p.180 in
Triebel, 1983). Since $\{U^{k}_{t}\} $ is also uniformly bounded
in $\mathcal{C}([0,T_{0}],B^{\sigma-1}_{2,1} )$(it
weak$^{\star}$-converges to $U_{t}$ in
$L^{\infty}([0,T_{0}],B^{\sigma-1}_{2,1})$), then $\{U^{k}\}$ is
uniformly bounded in Lip$([0,T_{0}],B^{\sigma-1}_{2,1})$, hence
uniformly equicontinuous on $[0,T_{0}]$ with the norm in $
B^{\sigma-1}_{2,1}$. By Proposition 2.2, Ascoli-Arzela theorem and
Cantor diagonal process, we deduce that
$$\phi U^{k}\rightarrow \phi U\ \ \ \mbox{in}\ \ \mathcal{C}([0,T_{0}],B^{\sigma-1}_{2,1}
)\ \ \mbox{as} \ \ k\rightarrow\infty, \ \ \ \mbox{for any}\ \
\phi\in\
 C_{c}^{\infty}.$$
The properties of strong convergence enable us to pass to the
limit in (4.10)-(4.11). Indeed, $U$ is a solution to (3.1)-(3.2).
Now, what remains is to check $U$ has the required regularity.
First, from $U\in \mathcal{C}([0,T_{0}],B^{\sigma-1}_{2,1} )$ and
an interpolation argument, we know  $U\in
\mathcal{C}([0,T_{0}],B^{\sigma'}_{2,1} )$ for any
$\sigma'<\sigma.$ Furthermore, we have $S_{q}U\in
\mathcal{C}([0,T_{0}],B^{\sigma}_{2,1})$ for any $q\in\mathbb{
N}.$ Then, combining with (4.18) (throw off the superscript $k$),
we derive that $S_{q}U$ converges uniformly to $U$ on $[0,T_{0}]$
with the norm in $B^{\sigma}_{2,1}$. This achieves $U\in
\mathcal{C}([0,T_{0}],B^{\sigma}_{2,1}).$ And then, using
Eq.(3.1), we see that $U_{t}\in
\mathcal{C}([0,T_{0}],B^{\sigma-1}_{2,1}),$ so $U(t,x)\in
\mathcal{C}^{1}([0,T_{0}]\times \mathbb{R}^{N})$. By virtue of the
flow map, we  get
$0<\theta'\leqslant\frac{\gamma-1}{2}m(t,x)+\bar{\psi}\leqslant\theta''$
for $(t,x)\in[0,T_{0}]\times \mathbb{R}^{N}$ according to Claim 1
($\theta',\theta'' $ are two positive constants).

(\textbf{Uniqueness}) \ Let $\tilde{m}=m_{1}-m_{2},\
\tilde{\textbf{u}}=\textbf{u}_{1}-\textbf{u}_{2},\
\tilde{\textbf{e}}=\textbf{e}_{1}-\textbf{e}_{2}$ where
$U_{1}=(m_{1},\textbf{u}_{1},\textbf{e}_{1})^{\top}$ and $
U_{2}=(m_{2},\textbf{u}_{2},\textbf{e}_{2})^{\top}$ are two
solutions to the system (3.1)-(3.2) with the same initial data
respectively. Then $\tilde{U}=(\tilde{m}, \tilde{\textbf{u}},
\tilde{\textbf{e}})^{\top}$ satisfies the following equations:
\begin{equation}\left\{
\begin{array}{l}\tilde{m}_{t}+\bar{\psi}\mbox{div}\tilde{\textbf{u}}=-\textbf{u}_{1}\cdot\nabla \tilde{m}-\tilde{\textbf{u}}\nabla
m_{2}-\frac{\gamma-1}{2}m_{1}\mbox{div}\tilde{\textbf{u}}-\frac{\gamma-1}{2}\tilde{m}\mbox{div}\textbf{u}_{2},\\
 \tilde{\textbf{u}}_{t}+\bar{\psi}\nabla \tilde{m}+\frac{\tilde{\textbf{u}}}{\tau}=-\textbf{u}_{1}\cdot\nabla\tilde{\textbf{u}}-\tilde{\textbf{u}}\nabla\textbf{u}_{2}-\frac{\gamma-1}{2}m_{1}\nabla\tilde{ m}-\frac{\gamma-1}{2}\tilde{ m}\nabla m_{2}+\tilde{\textbf{e}},
 \\
 \tilde{\textbf{e}}_{t}=-\nabla\Delta^{-1}\nabla\cdot\{\mathbf{H}(\tilde{m})\textbf{u}_{1}+(h(m_{2})+\bar{n})
 \tilde{\textbf{u}}\},
 \end{array}
\right.\end{equation} where the smooth function
$\mathbf{H}(\tilde{m}):=h(m_1)-h(m_2)$ satisfies
$\mathbf{H}(0)=0$.

 \noindent Similar to the derivation of (4.18), from Lemma 4.2 (take $p=2$ and $s=\sigma-1$), we obtain the
following estimate:
\begin{equation} \|\tilde{U}(t)\|_{B^{\sigma-1}_{2,1}}\leqslant C
\int^{t}_{0}\|\tilde{U}(\varsigma)\|_{B^{\sigma-1}_{2,1}}\Big(\|U_{1}(\varsigma)\|_{B^{\sigma}_{2,1}}+\|U_{2}(\varsigma)\|_{B^{\sigma}_{2,1}}\Big)d\varsigma\
, \ \ \ \mbox{for} \ \ t\in [0,T_{0}] . \end{equation} By
Gronwall's inequality, we have $\tilde{U}\equiv 0$.\end{proof}

\section{Global Existence and Exponential Stability}
In this section, we first state a proposition on the global
existence and exponential stability of classical solutions to
(3.1)-(3.2).
\begin{prop} Suppose  that  $ U_{0}\in B^{\sigma}_{2,1}$.
There  exists a  positive  constant
$\delta_{2}<\min\{\frac{\delta_{1}}{2C_{2}},\
\frac{\delta_{1}}{C_{1}}\}$ depending only on $A, \gamma, \tau$
and $\bar{n}$ such that if
$$\|U(\cdot,0)\|_{B^{\sigma}_{2,1}}\leqslant
\delta_{2}, $$ then there exists a unique global solution $U$ to
(3.1)-(3.2) satisfying
$$U\in \mathcal{C}([0,\infty),B^{\sigma}_{2,1})\cap \mathcal{C}^1([0,\infty),B^{\sigma-1}_{2,1})$$
and
$$\|U(\cdot,t)\|_{B^{\sigma}_{2,1}}+\|U_{t}(\cdot,t)\|_{B^{\sigma-1}_{2,1}}\leqslant
C_{1}\|U(\cdot,0)\|_{B^{\sigma}_{2,1}}\exp(-\mu_{1}t),\ \ \
t\geqslant 0,$$ where $\delta_{1},\ \mu_{1}$, $C_{1}$ are three
positive constants given by Proposition 5.2 and the positive
constant $C_{2}$ is given in (5.19),
$U=(m,\textbf{u},\textbf{e})^{\top}$ and
$U_{t}=(m_{t},\textbf{u}_{t},\textbf{e}_{t})^{\top}$.
\end{prop}

The proof of above proposition mainly depends on a crucial
\textit{a-priori} estimate (Proposition 5.2). To do this, we need
some lemmas.
\begin{lem} If $(m,\textbf{u},\textbf{e})\in
\mathcal{C}([0,T],B^{\sigma}_{2,1})\cap
\mathcal{C}^1([0,T],B^{\sigma-1}_{2,1})$ is a solution of Eq.(3.1)
for any given $T>0$, then
\begin{eqnarray}&&\frac{d}{dt}\Big(\|\Delta_{q}m_{t}\|^2_{L^2}+\|\Delta_{q}\textbf{u}_{t}\|^2_{L^2}+\frac{1}{\bar{n}}\|\Delta_{q}\textbf{e}_{t}\|^2_{L^2}\Big)+\frac{2}{\tau}\|\Delta_{q}\textbf{u}_{t}\|^2_{L^2}\nonumber\\&\leqslant&
2\|\textbf{u}_{t}\|_{L^{\infty}}(\|\Delta_{q}\nabla
m\|_{L^2}\|\Delta_{q}m_{t}\|_{L^2}+\|\Delta_{q}\nabla
\textbf{u}\|_{L^2}\|\Delta_{q}\textbf{u}_{t}\|_{L^2})+\|\nabla\textbf{u}\|_{L^{\infty}}(\|\Delta_{q}m_{t}\|^2_{L^2}\nonumber\\
&&+\|\Delta_{q}\textbf{u}_{t}\|^2_{L^2})+(\gamma-1)\|m_{t}\|_{L^{\infty}}(\|\Delta_{q}\mathrm{div}\textbf{u}\|_{L^2}
\|\Delta_{q}m_{t}\|_{L^2}+\|\Delta_{q}\nabla
m\|_{L^2}\|\Delta_{q}\textbf{u}_{t}\|_{L^2})\nonumber\\
&&+(\gamma-1)\|\nabla
m\|_{L^{\infty}}\|\Delta_{q}m_{t}\|_{L^2}\|\Delta_{q}\textbf{u}_{t}\|_{L^2}+(\|[\textbf{u}_{t},
\Delta_{q}]\cdot\nabla
m\|_{L^2}+\|[\textbf{u},\Delta_{q}]\cdot\nabla
m_{t}\|_{L^2})\nonumber\\
&&\times\|\Delta_{q}m_{t}\|_{L^2}+(\|[\textbf{u}_{t},
\Delta_{q}]\cdot\nabla \textbf{u}\|_{L^2}+\|[\textbf{u},
\Delta_{q}]\cdot\nabla
\textbf{u}_{t}\|_{L^2})\|\Delta_{q}\textbf{u}_{t}\|_{L^2}+(\gamma-1)\nonumber\\
&&\times(\|[m_{t},\Delta_{q}]\mathrm{div}\textbf{u}\|_{L^2}+\|[m,\Delta_{q}]\mathrm{div}\textbf{u}_{t}\|_{L^2})\|\Delta_{q}m_{t}\|_{L^2}+(\gamma-1)(\|[m_{t},\Delta_{q}]\nabla
m\|_{L^2}\nonumber\\
&&+\|[m,\Delta_{q}]\nabla
m_{t}\|_{L^2})\|\Delta_{q}\textbf{u}_{t}\|_{L^2}
+\frac{2}{\bar{n}}\Big(\|\Delta_{q}(\mathcal{H}(m)m_{t}\textbf{u})\|_{L^2}\nonumber\\&&+(A\gamma)^{-\frac{1}{2}}\bar{n}^{
 \frac{3-\gamma}{2}}\|\Delta_{q}(m_{t}\textbf{u})\|_{L^2}+\|\Delta_{q}(h(m)\textbf{u}_{t})\|_{L^2}\Big)\|\Delta_{q}\textbf{e}_{t}\|_{L^2},
\end{eqnarray}
where
$\mathcal{H}(m)=(A\gamma)^{-\frac{1}{2}}\{(A\gamma)^{-\frac{1}{2}}(\frac{\gamma-1}{2}m+\bar{\psi})\}^{
 \frac{3-\gamma}{\gamma-1}}-(A\gamma)^{-\frac{1}{2}}\bar{n}^{
 \frac{3-\gamma}{2}}$ is a smooth function on $\{m|\frac{\gamma-1}{2}m+\bar{\psi}>0\}$ satisfying $\mathcal{H}(0)=0$.
\end{lem}
\begin{proof} By differentiating the first two equations of
Eq.(4.1) with respect to variable $t$ once, integrating them over
$\mathbb{R}^{N}$ after multiplying
$\Delta_{q}m_{t},\Delta_{q}\textbf{u}_{t}$ respectively, similar
to the derivation of (4.16), through tedious but straightforward
calculations, we can obtain (5.1).
\end{proof}
In addition, we give some auxiliary estimates, which are divided
into high- and low- frequency cases.
\begin{lem}If $(m,\textbf{u},\textbf{e})\in\mathcal{C}([0,T],B^{\sigma}_{2,1})\cap
 \mathcal{C}^1([0,T],B^{\sigma-1}_{2,1})$ is a solution of Eq.(3.1)
for any given $T>0$, then
\begin{eqnarray}
&&\|\Delta_{q}m_{t}\|^2_{L^2}\nonumber\\
&\leqslant&
\Big(\bar{C}\bar{\psi}2^{q}\|\Delta_{q}\textbf{u}\|_{L^2}+\|\textbf{u}\|_{L^{\infty}}\|\Delta_{q}\nabla
m\|_{L^2}+\frac{\gamma-1}{2}\|m\|_{L^{\infty}}\|\Delta_{q}\mathrm{div}
\textbf{u}\|_{L^2}\nonumber\\&&+\|[\textbf{u},\Delta_{q}]\nabla
m\|_{L^2}+\frac{\gamma-1}{2}\|[m,\Delta_{q}]\mathrm{div}
\textbf{u}\|_{L^2}\Big)\|\Delta_{q}m_{t}\|_{L^2} \ (q\geqslant-1);
\end{eqnarray}
\begin{eqnarray}&&2^{q}\bar{\psi}\|\Delta_{q}m\|^2_{L^2}\nonumber\\&\leqslant& \bar{C}\Big(\frac{1}{\tau}\|\Delta_{q}\textbf{u}\|_{L^2}+\|\Delta_{q}\textbf{u}_{t}\|_{L^2}+\|\textbf{u}\|_{L^{\infty}}\|\Delta_{q}\nabla\textbf{u}\|_{L^2}\nonumber\\&&+\|[\textbf{u},\Delta_{q}]\cdot\nabla \textbf{u}\|_{L^2}+\frac{\gamma-1}{2}
\|m\|_{L^{\infty}}\|\Delta_{q}\nabla m\|_{L^2}\nonumber\\&&+
\frac{\gamma-1}{2}\|[m,\Delta_{q}]\cdot\nabla
m\|_{L^2}+\|\Delta_{q}(\tilde{h}(m)m)\|_{L^2}\Big)\|\Delta_{q}m\|_{L^2}\
 (q\geqslant0) ;\end{eqnarray}
\begin{eqnarray}&&(A\gamma)^{-\frac{1}{2}}\bar{n}^{\frac{3-\gamma}{2}} \|\Delta_{-1}m\|^2_{L^2}\nonumber\\ &\leqslant& \bar{C}\Big(\frac{1}{\tau}\|\Delta_{-1}\textbf{u}\|_{L^2}+\|\Delta_{-1}\textbf{u}_{t}\|_{L^2}+\|\textbf{u}\|_{L^{\infty}}\|\Delta_{-1}\nabla\textbf{u}\|_{L^2}
\nonumber\\&&+\|[\textbf{u},\Delta_{-1}]\cdot\nabla
\textbf{u}\|_{L^2} +\frac{\gamma-1}{2}
\|m\|_{L^{\infty}}\|\Delta_{-1}\nabla
m\|_{L^2}\nonumber\\&&+\frac{\gamma-1}{2}\|[m,\Delta_{-1}]\cdot\nabla
m\|_{L^2}+\|\Delta_{-1}(\tilde{h}(m)m)\|_{L^2}\Big)\|\Delta_{-1}m\|_{L^2}
;\end{eqnarray}
\begin{eqnarray}&&2^{q}\|\Delta_{q}\textbf{e}\|^2_{L^2}\nonumber\\ &\leqslant& \bar{C}\{(A\gamma)^{-\frac{1}{2}}\bar{n}^{\frac{3-\gamma}{2}}\|\Delta_{q}m\|_{L^2}+\|\Delta_{q}(\tilde{h}(m)m)\|_{L^2}\}\|\Delta_{q}\textbf{e}\|_{L^2}\ (q\geqslant0);\end{eqnarray}
\begin{eqnarray}&&\|\Delta_{-1}\textbf{e}\|^2_{L^2}\nonumber\\&\leqslant&
 \Big(\|\Delta_{-1}\textbf{u}_{t}\|_{L^2}+\frac{1}{\tau}\|\Delta_{-1}\textbf{u}\|_{L^2}+\bar{C}\bar{\psi}\|\Delta_{-1}
 m\|_{L^2}+\|\Delta_{-1}(\textbf{u}\cdot\nabla
 \textbf{u})\|_{L^2}\nonumber\\&&+\frac{\gamma-1}{2}\|\Delta_{-1}(m\nabla
 m)\|_{L^2}\Big)\|\Delta_{-1}\textbf{e}\|_{L^2};\end{eqnarray}
\begin{equation}\ \ \|\Delta_{q}\textbf{e}_{t}\|^2_{L^2}\leqslant
\bar{C}(\|\Delta_{q}(h(m)\textbf{u})\|_{L^2}+\bar{n}\|\Delta_{q}\textbf{u}\|_{L^2})\|\Delta_{q}\textbf{e}_{t}\|_{L^2}\
\ (q\geqslant-1),\end{equation} where
$\tilde{h}(m)=\int_{0}^{1}h'(\varsigma
m)d\varsigma-(A\gamma)^{-\frac{1}{2}}\bar{n}^{\frac{3-\gamma}{2}}$
is a smooth function on $\{m|\frac{\gamma-1}{2}\varsigma
m+\bar{\psi}>0,\ \varsigma\in[0,1]\}$, $\tilde{h}(0)=0$ and
$\bar{C}$ is a uniform constant independent of $A, \gamma,\tau$
and $\bar{n}$.
\end{lem}
\begin{proof} (1) Using the first equation of Eq.(3.1), we have
\begin{equation}m_{t}=-\Big(\bar{\psi}\mbox{div}\textbf{u}+\textbf{u}\cdot\nabla
m+\frac{\gamma-1}{2}m\mbox{div}\textbf{u}\Big).\end{equation} By
applying the operator $\Delta_{q}(q\geqslant-1)$ to (5.8),
integrating it over $\mathbb{R}^{N}$ after multiplying
 $\Delta_{q}m_{t}$, we can arrive at (5.2) with the aid of H\"{o}lder's
inequality.

(2) Using the second equation of Eq.(3.1), we get
\begin{equation}\bar{\psi}\nabla
m=-\Big(\textbf{u}_{t}+\frac{\textbf{u}}{\tau}+\textbf{u}\cdot\nabla\textbf{u}+\frac{\gamma-1}{2}m\nabla
m-\textbf{e}\Big).\end{equation} By Lemma 4.1, we have
$\mbox{div}\textbf{e}=h(m)$. Integration by parts gives
\begin{eqnarray*}\int\Delta_{q}\textbf{e}\cdot\nabla\Delta_{q}m&=&-\int\Delta_{q}\mbox{div}\textbf{e}\Delta_{q}m=-\int
\Delta_{q}h(m)\Delta_{q}m\\&=&-\int\Delta_{q}(h(m)-h(0))\Delta_{q}m\\&=&-(A\gamma)^{-\frac{1}{2}}\bar{n}^{\frac{3-\gamma}{2}}\|\Delta_{q}m\|^2_{L^2}-\int\Delta_{q}(\tilde{h}(m)m)\Delta_{q}m,
\end{eqnarray*}
where  $\tilde{h}(m)=\int_{0}^{1}h'(\varsigma
m)d\varsigma-(A\gamma)^{-\frac{1}{2}}\bar{n}^{\frac{3-\gamma}{2}}$
is a smooth function on $\{m|\frac{\gamma-1}{2}\varsigma
m+\bar{\psi}>0,\ \varsigma\in[0,1]\}$ satisfying $
\tilde{h}(0)=0$. By applying the operator $\Delta_{q}$ to (5.9),
integrating it over $\mathbb{R}^{N}$ after multiplying
$\Delta_{q}\nabla m$, we can obtain
\begin{eqnarray}&&\bar{\psi}\|\Delta_{q}\nabla m\|^2_{L^2}+(A\gamma)^{-\frac{1}{2}}\bar{n}^{\frac{3-\gamma}{2}}\|\Delta_{q}
m\|^2_{L^2}\nonumber\\&\leqslant&
\Big(\frac{1}{\tau}\|\Delta_{q}\textbf{u}\|_{L^2}+\|\Delta_{q}\textbf{u}_{t}\|_{L^2}+\|\textbf{u}\|_{L^{\infty}}\|\Delta_{q}\nabla\textbf{u}\|_{L^2}+\|[\textbf{u},\Delta_{q}]\cdot\nabla
\textbf{u}\|_{L^2}\nonumber\\&&+\frac{\gamma-1}{2}
\|\textbf{m}\|_{L^{\infty}}\|\Delta_{q}\nabla
m\|_{L^2}+\frac{\gamma-1}{2}\|[\textbf{m},\Delta_{q}]\cdot\nabla
m\|_{L^2}\Big)\|\nabla\Delta_{q}m\|_{L^2}\nonumber\\&&+\|\Delta_{q}(\tilde{h}(m)m)\|_{L^2}\|\Delta_{q}m\|_{L^2}.\end{eqnarray}
By Lemma 2.1, we know
$$\|\Delta_{q}\nabla m\|_{L^2}\approx
2^{q}\|\Delta_{q}m\|_{L^2}\ \ \ (q\geqslant0),$$ so (5.3) follows.
For the low frequency $(q=-1)$, by (5.10), we get (5.4)
immediately.

(3) By applying the operator $\Delta_{q}$ to both sides of
$\mbox{div}\textbf{e}=h(m)$ ($q\geqslant0$), integrating it over $
\mathbb{R}^{N}$ after multiplying
$\Delta_{q}\mbox{div}\textbf{e}$, we
   reach (5.5) in virtue of
H\"{o}lder's inequality.  For $q=-1$, using the second equation of
Eq.(3.1), we can obtain (5.6).

(4) From the $L^2$- boundedness of Riesz transform, we can derive
(5.7) directly from the last equation of Eq.(3.1).\end{proof} The
crucial \textit{a-priori} estimate is comprised in the following
proposition.
\begin{prop}
There exist three positive constants $\delta_{1}, C_{1}$ and
$\mu_{1}$ depending only on $A, \gamma, \bar{n}$ and $\tau$ such
that for any $T>0$, if
\begin{eqnarray}\sup_{0\leqslant t\leqslant
T}\Big(\|U(\cdot,t)\|_{B^{\sigma}_{2,1}}+\|U_{t}(\cdot,t)\|_{B^{\sigma-1}_{2,1}}\Big)\leqslant
\delta_{1},\end{eqnarray} then
\begin{eqnarray}\|U(\cdot,t)\|_{B^{\sigma}_{2,1}}+\|U_{t}(\cdot,t)\|_{B^{\sigma-1}_{2,1}}\leqslant
C_{1}\|U(\cdot,0)\|_{B^{\sigma}_{2,1}}\exp(-\mu_{1}t), \ \
t\in[0,T].\nonumber
\end{eqnarray}
\end{prop}
\begin{proof} From the \textit{a-priori} assumption (5.11), we
have
\begin{eqnarray}\sup_{0\leqslant t\leqslant T} \Big(\|U(\cdot,t)||_{W^{1,\infty}}+\|U_{t}(\cdot,t)\|_{L^{\infty}}\Big)\leqslant
C\delta_{1}. \end{eqnarray} To ensure the smoothness of functions
$h(m),\mathcal{H}(m)$ and $\tilde{h}(m)$, we may choose\\
$0<\delta_{1}\leqslant\frac{\bar{\psi}}{(\gamma-1)C}$, then
$$\frac{\gamma-1}{2}m(t,x)+\bar{\psi}\geqslant \frac{\bar{\psi}}{2}>0,\ \ \ (t,x)\in[0,T]\times\mathbb{R}^{N}$$
and
$$\frac{\gamma-1}{2}\varsigma m(t,x)+\bar{\psi}\geqslant\frac{\bar{\psi}}{2}>0, \ \varsigma\in[0,1],\ (t,x)\in[0,T]\times\mathbb{R}^{N}.$$
 From (4.2)-(4.4), we set
\begin{eqnarray}
I_{1,q}&:=&\|\nabla\textbf{u}\|_{L^{\infty}}(||\Delta_{q}m\|^2_{L^2}+||\Delta_{q}\textbf{u}\|^2_{L^2})
+2\|[\textbf{u},\Delta_{q}]\cdot\nabla
m\|_{L^2}||\Delta_{q}m\|_{L^2}\nonumber\\&&\hspace{-2.5mm}+2\|[\textbf{u},\Delta_{q}]\cdot\nabla
\textbf{u}\|_{L^2}||\Delta_{q}\textbf{u}\|_{L^2}+(\gamma-1)\|\nabla
m\|_{L^{\infty}}||\Delta_{q}m\|_{L^2}||\Delta_{q}\textbf{u}\|_{L^2}\nonumber\\&&\hspace{-2.5mm}+(\gamma-1)\|[m,\Delta_{q}]\nabla
m\|_{L^2}||\Delta_{q}\textbf{u}\|_{L^2}
+(\gamma-1)\|[m,\Delta_{q}]\mathrm{div}\textbf{u}\|_{L^2}||\Delta_{q}m\|_{L^2}\nonumber\\&&\hspace{-2.5mm}+\frac{2}{\bar{n}}\|\Delta_{q}(h(m)\textbf{u})\|_{L^2}\|\Delta_{q}\textbf{e}\|_{L^2}.
\end{eqnarray}
$I_{2,q}$ denotes the right-hand side of inequality (5.1). For the
proof of Proposition 5.2, we are going to divided it into the
following two lemmas.
\begin{lem}$(q\geqslant0)$
There exists a positive constant $\mu_{2}$ depending only on $A,
\gamma, \tau$ and $\bar{n}$ such that the following estimate
holds:
\begin{eqnarray}
&&2^{q(\sigma-1)}\frac{d}{dt}\Big\{2^{2q}\Big(\|\Delta_{q}m\|^2_{L^2}+\|\Delta_{q}\textbf{u}\|^2_{L^2}+\frac{1}{\bar{n}}\|\Delta_{q}\textbf{e}\|^2_{L^2}\Big)+\Big(\|\Delta_{q}m_{t}\|^2_{L^2}\nonumber\\&&+\|\Delta_{q}\textbf{u}_{t}\|^2_{L^2}+\frac{1}{\bar{n}}\|\Delta_{q}\textbf{e}_{t}\|^2_{L^2}\Big)\Big\}^{\frac{1}{2}}
+\mu_{2}(2^{q\sigma}\|\Delta_{q}U\|_{L^2}+2^{q(\sigma-1)}\|\Delta_{q}U_{t}\|_{L^2})\nonumber\\
&\leqslant& C\{2^{q(\sigma-1)}(\|
U\|_{W^{1,\infty}}+\|U_{t}\|_{L^{\infty}})(2^{q}\|\Delta_{q}U\|_{L^2}+\|\Delta_{q}U_{t}\|_{L^2})
+c_{q}(\|U\|_{B^{\sigma}_{2,1}}\nonumber\\&&+\|U_{t}\|_{B^{\sigma-1}_{2,1}})\|U\|_{B^{\sigma}_{2,1}}
+2^{q\sigma}\|\Delta_{q}(h(m)\textbf{u})\|_{L^2}
+2^{q(\sigma-1)}(\|\Delta_{q}(\mathcal{H}(m)m_{t}\textbf{u})\|_{L^2}\nonumber\\&&+\|\Delta_{q}(m_{t}\textbf{u})\|_{L^2}+\|\Delta_{q}(h(m)\textbf{u}_{t})\|_{L^2}+\|\Delta_{q}(\tilde{h}(m)m)\|_{L^2})\},
\end{eqnarray}
where $C>0$ is a harmless constant depending only on $A, \gamma,
\tau$ and $\bar{n}$.
\end{lem}
\noindent\textbf{Proof of Lemma 5.3}. Combining (4.2)-(4.4), Lemma
5.1 and 5.2, we have $(q\geqslant0)$
\begin{eqnarray}
&&\hspace{-3mm}\frac{d}{dt}\Big\{2^{2q}\Big(\|\Delta_{q}m\|^2_{L^2}+\|\Delta_{q}\textbf{u}\|^2_{L^2}+\frac{1}{\bar{n}}\|\Delta_{q}\textbf{e}\|^2_{L^2}\Big)
+\Big(\|\Delta_{q}m_{t}\|^2_{L^2}+\|\Delta_{q}\textbf{u}_{t}\|^2_{L^2}\nonumber\\&&+\frac{1}{\bar{n}}\|\Delta_{q}\textbf{e}_{t}\|^2_{L^2}\Big)\Big\}
+\beta_{1}\bar{\psi}2^{2q}\|\Delta_{q}m\|^2_{L^2}+\frac{2}{\tau}2^{2q}\|\Delta_{q}\textbf{u}\|^2_{L^2}+\frac{\beta_{2}}{\bar{n}}2^{2q}\|\Delta_{q}\textbf{e}\|^2_{L^2}
\nonumber\\&&+\beta_{3}\|\Delta_{q}m_{t}\|^2_{L^2}+\frac{2}{\tau}\|\Delta_{q}\textbf{u}_{t}\|^2_{L^2}+\frac{\beta_{4}}{\bar{n}}\|\Delta_{q}\textbf{e}_{t}\|^2_{L^2}
\nonumber\\&\leqslant& 2^{2q}I_{1,q}
+I_{2,q}+\beta_{1}\bar{C}2^{q}\Big(\frac{1}{\tau}\|\Delta_{q}\textbf{u}\|_{L^2}+\|\Delta_{q}\textbf{u}_{t}\|_{L^2}+\|\textbf{u}\|_{L^{\infty}}\|\Delta_{q}\nabla
\textbf{u}\|_{L^2}\nonumber\\&&+\|[\textbf{u},\Delta_{q}]\nabla
\textbf{u}\|_{L^2}+\frac{\gamma-1}{2}\|m\|_{L^{\infty}}\|\Delta_{q}\nabla
m\|_{L^2}+\frac{\gamma-1}{2}\|[m,\Delta_{q}]\nabla
m\|_{L^2}\nonumber\\&&+\|\Delta_{q}(\tilde{h}(m)m)\|_{L^2}\Big)\|\Delta_{q}m\|_{L^2}+\frac{\beta_{2}\bar{C}}{\bar{n}}2^{q}\Big((A\gamma)^{-\frac{1}{2}}\bar{n}^{\frac{3-\gamma}{2}}\|\Delta_{q}m\|_{L^2}\nonumber\\&&+\|\Delta_{q}(\tilde{h}(m)m)\|_{L^2}\Big)\|\Delta_{q}\textbf{e}\|_{L^2}
+\beta_{3}
\Big(\bar{C}\bar{\psi}2^{q}\|\Delta_{q}\textbf{u}\|_{L^2}+\|\textbf{u}\|_{L^{\infty}}\|\Delta_{q}\nabla
m\|_{L^2}\nonumber\\&&+\frac{\gamma-1}{2}\|m\|_{L^{\infty}}\|\Delta_{q}\mathrm{div}
\textbf{u}\|_{L^2}+\frac{\gamma-1}{2}\|[m,\Delta_{q}]\mathrm{div}
\textbf{u}\|_{L^2}+\|[\textbf{u},\Delta_{q}]\nabla
m\|_{L^2}\Big)\nonumber\\&&\times\|\Delta_{q}m_{t}\|_{L^2}+\frac{\beta_{4}\bar{C}}{\bar{n}}\Big(\|\Delta_{q}h(m)\textbf{u}\|_{L^2}+\bar{n}\|\Delta_{q}\textbf{u}\|_{L^2}\Big)\|\Delta_{q}\textbf{e}_{t}\|_{L^2},
\end{eqnarray}
where these positive constants $\beta_{1}, \beta_{2}, \beta_{3}$
and $\beta_{4}$ satisfy
$$\beta_{1}=\min\Big\{\frac{\bar{\psi}\tau}{2\bar{C}^2},\ \frac{\bar{\psi}}{\tau\bar{C}^2}\Big\},\ \ \beta_{2}=\frac{A\gamma\bar{\psi}}{2\bar{C}^2\bar{n}^{2-\gamma}}\beta_{1},\ \
\beta_{3}=\frac{1}{\tau\bar{C}^2\bar{\psi}^2},\ \
\beta_{4}=\frac{1}{\tau\bar{n}\bar{C}^2},$$respectively. We
introduce them in order to eliminate quadratic terms in the
right-hand side of (5.15). There are no quadratic terms in
$I_{1,q}$ and $I_{2,q}$. By Young's inequality, the first
quadratic term can be estimated in this way:
\begin{eqnarray*}
&&\beta_{1}\bar{C}\frac{1}{\tau}2^{q}\|\Delta_{q}\textbf{u}\|_{L^2}\|\Delta_{q}m\|_{L^2}\\
&\leqslant&\sqrt{\frac{\beta_{1}}{\bar{\psi}}}\bar{C}\frac{1}{\tau}2^{q}\|\Delta_{q}\textbf{u}\|_{L^2}\cdot
2^{q}\sqrt{\beta_{1}\bar{\psi}}\|\Delta_{q}m\|_{L^2}\\
&\leqslant&\frac{\bar{C}^2\beta_{1}}{\bar{\psi}\tau^2}2^{2q}\|\Delta_{q}\textbf{u}\|^2_{L^2}+\frac{1}{4}\beta_{1}\bar{\psi}2^{2q}\|\Delta_{q}m\|^2_{L^2}\\
&\leqslant&\frac{1}{2\tau}2^{2q}\|\Delta_{q}\textbf{u}\|^2_{L^2}+\frac{1}{4}\beta_{1}\bar{\psi}2^{2q}\|\Delta_{q}m\|^2_{L^2}.\end{eqnarray*}
The remainder quadratic terms in the right-hand side of (5.15) are
estimated similarly as follows:
\begin{eqnarray*}&&\beta_{1}\bar{C}2^{q}\|\Delta_{q}\textbf{u}_{t}\|_{L^2}\|\Delta_{q}m\|_{L^2}\\&\leqslant&\frac{\beta_{1}\bar{C}^2}{\bar{\psi}}\|\Delta_{q}\textbf{u}_{t}\|^2_{L^2}+\frac{1}{4}\beta_{1}\bar{\psi}2^{2q}\|\Delta_{q}m\|^2_{L^2}\\
&\leqslant&\frac{1}{\tau}\|\Delta_{q}\textbf{u}_{t}\|^2_{L^2}+\frac{1}{4}\beta_{1}\bar{\psi}2^{2q}\|\Delta_{q}m\|^2_{L^2};
\end{eqnarray*}
\begin{eqnarray*}&&\frac{\beta_{2}\bar{C}}{\bar{n}}(A\gamma)^{-\frac{1}{2}}\bar{n}^{\frac{3-\gamma}{2}}2^{q}\|\Delta_{q}m\|_{L^2}\|\Delta_{q}\textbf{e}\|_{L^2}\\
&\leqslant&\frac{\bar{n}^{3-\gamma}\beta_{2}\bar{C}^2}{2A\gamma\bar{n}}2^{2q}\|\Delta_{q}m\|^2_{L^2}+\frac{\beta_{2}}{2\bar{n}}2^{2q}\|\Delta_{q}\textbf{e}\|^2_{L^2}\\
&=&\frac{1}{4}\beta_{1}\bar{\psi}2^{2q}\|\Delta_{q}m\|^2_{L^2}+\frac{\beta_{2}}{2\bar{n}}2^{2q}\|\Delta_{q}\textbf{e}\|^2_{L^2};
\end{eqnarray*}
\begin{eqnarray*}&&\bar{C}\bar{\psi}\beta_{3}2^{q}\|\Delta_{q}\textbf{u}\|_{L^2}\|\Delta_{q}m_{t}\|_{L^2}\\&\leqslant&
\frac{\bar{C}^2\bar{\psi}^2}{2}\beta_{3}2^{2q}\|\Delta_{q}\textbf{u}\|^2_{L^2}+\frac{1}{2}\beta_{3}\|\Delta_{q}m_{t}\|^2_{L^2}\\
&=&\frac{1}{2\tau}2^{2q}\|\Delta_{q}\textbf{u}\|^2_{L^2}+\frac{1}{2}\beta_{3}\|\Delta_{q}m_{t}\|^2_{L^2};
\end{eqnarray*}
\begin{eqnarray*}
&&\beta_{4}\bar{C}\|\Delta_{q}\textbf{u}\|_{L^2}\|\Delta_{q}\textbf{e}_{t}\|_{L^2}\\&\leqslant&\frac{\bar{n}\bar{C}^2}{2}\beta_{4}2^{2q}\|\Delta_{q}\textbf{u}\|^2_{L^2}+\frac{1}{2\bar{n}}\beta_{4}\|\Delta_{q}\textbf{e}_{t}\|^2_{L^2}\\
&=&\frac{1}{2\tau}2^{2q}\|\Delta_{q}\textbf{u}\|^2_{L^2}+\frac{1}{2\bar{n}}\beta_{4}\|\Delta_{q}\textbf{e}_{t}\|^2_{L^2}.
\end{eqnarray*}
Then (5.15) becomes into
\begin{eqnarray}
&&\frac{d}{dt}\Big\{2^{2q}\Big(\|\Delta_{q}m\|^2_{L^2}+\|\Delta_{q}\textbf{u}\|^2_{L^2}+\frac{1}{\bar{n}}\|\Delta_{q}\textbf{e}\|^2_{L^2}\Big)+\Big(\|\Delta_{q}m_{t}\|^2_{L^2}+\|\Delta_{q}\textbf{u}_{t}\|^2_{L^2}\nonumber\\&&+\frac{1}{\bar{n}}\|\Delta_{q}\textbf{e}_{t}\|^2_{L^2}\Big)\Big\}
+\frac{\beta_{1}}{4}\bar{\psi}2^{2q}\|\Delta_{q}m\|^2_{L^2}+\frac{1}{2\tau}2^{2q}\|\Delta_{q}\textbf{u}\|^2_{L^2}+\frac{\beta_{2}}{2\bar{n}}2^{2q}\|\Delta_{q}\textbf{e}\|^2_{L^2}
\nonumber\\&&+\frac{\beta_{3}}{2}\|\Delta_{q}m_{t}\|^2_{L^2}+\frac{1}{\tau}\|\Delta_{q}\textbf{u}_{t}\|^2_{L^2}+\frac{\beta_{4}}{2\bar{n}}\|\Delta_{q}\textbf{e}_{t}\|^2_{L^2}
\nonumber\\&\leqslant& 2^{2q}I_{1,q}
+I_{2,q}+\beta_{1}\bar{C}2^{q}\Big(\|\textbf{u}\|_{L^{\infty}}\|\Delta_{q}\nabla
\textbf{u}\|_{L^2}+\frac{\gamma-1}{2}\|m\|_{L^{\infty}}\|\Delta_{q}\nabla
m\|_{L^2}\nonumber\\&&+\|[\textbf{u},\Delta_{q}]\nabla
\textbf{u}\|_{L^2}+\frac{\gamma-1}{2}\|[m,\Delta_{q}]\nabla
m\|_{L^2}+\|\Delta_{q}(\tilde{h}(m)m)\|_{L^2}\Big)\|\Delta_{q}m\|_{L^2}\nonumber\\&&+
\frac{\beta_{2}\bar{C}}{\bar{n}}2^{q}\|\Delta_{q}(\tilde{h}(m)m)\|_{L^2}\|\Delta_{q}\textbf{e}\|_{L^2}+\beta_{3}
\Big(\|\textbf{u}\|_{L^{\infty}}\|\Delta_{q}\nabla
m\|_{L^2}+\frac{\gamma-1}{2}\|m\|_{L^{\infty}}\nonumber\\&&\times\|\Delta_{q}\mathrm{div}
\textbf{u}\|_{L^2}+\|[\textbf{u},\Delta_{q}]\nabla
m\|_{L^2}+\frac{\gamma-1}{2}\|[m,\Delta_{q}]\mathrm{div}
\textbf{u}\|_{L^2}\Big)\|\Delta_{q}m_{t}\|_{L^2}\nonumber\\&&+\frac{\beta_{4}\bar{C}}{\bar{n}}\|\Delta_{q}h(m)\textbf{u}\|_{L^2}\|\Delta_{q}\textbf{e}_{t}\|_{L^2}.
\end{eqnarray}
Dividing (5.16) by
$$\Big\{2^{2q}\Big(\|\Delta_{q}m\|^2_{L^2}+\|\Delta_{q}\textbf{u}\|^2_{L^2}+\frac{1}{\bar{n}}\|\Delta_{q}\textbf{e}\|^2_{L^2}\Big)+\Big(\|\Delta_{q}m_{t}\|^2_{L^2}+\|\Delta_{q}\textbf{u}_{t}\|^2_{L^2}+\frac{1}{\bar{n}}\|\Delta_{q}\textbf{e}_{t}\|^2_{L^2}\Big)\Big\}^{\frac{1}{2}}$$
and multiplying (5.16) by the factor $2^{q(\sigma-1)}$, we get
(5.14) immediately with the help of Lemma 4.2, which completes the
proof of Lemma 5.3. \ \ \hspace{27mm}$\square$

For the case of low frequency ($q=-1$), we also have the following
\textit{a-priori} estimate in a similar way.
\begin{lem}$(q=-1)$
There exists a positive constant $\mu_{3}$ depending only on $A,
\gamma, \tau$ and $\bar{n}$ such that the following estimate
holds:
\begin{eqnarray}
&&\frac{1}{2^{\sigma-1}}\frac{d}{dt}\Big\{\frac{1}{4}\Big(\|\Delta_{-1}m\|^2_{L^2}+\|\Delta_{-1}\textbf{u}\|^2_{L^2}+\frac{1}{\bar{n}}\|\Delta_{-1}\textbf{e}\|^2_{L^2}\Big)+\Big(\|\Delta_{-1}m_{t}\|^2_{L^2}
\nonumber\\&&+\|\Delta_{-1}\textbf{u}_{t}\|^2_{L^2}+\frac{1}{\bar{n}}\|\Delta_{-1}\textbf{e}_{t}\|^2_{L^2}\Big)\Big\}^{\frac{1}{2}}
+\mu_{2}\Big(\frac{1}{2^{\sigma}}\|\Delta_{-1}U\|_{L^2}+\frac{1}{2^{\sigma-1}}\|\Delta_{-1}U_{t}\|_{L^2}\Big)\nonumber\\
&\leqslant& C\Big\{\frac{1}{2^{\sigma-1}}\Big(\|
U\|_{W^{1,\infty}}+\|U_{t}\|_{L^{\infty}}\Big)\Big(\frac{1}{2}\|\Delta_{-1}U\|_{L^2}+\|\Delta_{-1}U_{t}\|_{L^2}\Big)
+c_{-1}(\|U\|_{B^{\sigma}_{2,1}}\nonumber\\&&+\|U_{t}\|_{B^{\sigma-1}_{2,1}})\|U\|_{B^{\sigma}_{2,1}}
+\frac{1}{2^{\sigma}}\|\Delta_{-1}(h(m)\textbf{u})\|_{L^2}
+\frac{1}{2^{\sigma-1}}\Big(\|\Delta_{-1}(\mathcal{H}(m)m_{t}\textbf{u})\|_{L^2}\nonumber\\&&+\|\Delta_{-1}(m_{t}\textbf{u})\|_{L^2}+\|\Delta_{-1}(h(m)\textbf{u}_{t})\|_{L^2}+\|\Delta_{-1}(\tilde{h}(m)m)\|_{L^2}\Big)\Big\},
\end{eqnarray}
where $C>0$ is a harmless constant depending only on $A, \gamma,
\tau$ and $\bar{n}$.
\end{lem}
 Summing (5.14) on $q \in \mathbb{N}\cup\{0\}$ and adding (5.17) together, according to Proposition 2.3 and \textit{a-priori} assumption
(5.11)-(5.12), we get the following differential inequality:
\begin{eqnarray}\frac{d}{dt}Q+\mu_{4}\Big(\|U\|_{B^{\sigma}_{2,1}}+\|U_{t}\|_{B^{\sigma-1}_{2,1}}\Big)\leqslant
C\delta_{1}\Big(\|U\|_{B^{\sigma}_{2,1}}+\|U_{t}\|_{B^{\sigma-1}_{2,1}}\Big),\end{eqnarray}
where
\begin{eqnarray}Q:&=&\sum_{q\geqslant-1}2^{q(\sigma-1)}\Big\{2^{2q}\Big(\|\Delta_{q}m\|^2_{L^2}+\|\Delta_{q}\textbf{u}\|^2_{L^2}+\frac{1}{\bar{n}}\|\Delta_{q}\textbf{e}\|^2_{L^2}\Big)\nonumber\\&&+\Big(\|\Delta_{q}m_{t}\|^2_{L^2}+\|\Delta_{q}\textbf{u}_{t}\|^2_{L^2}+\frac{1}{\bar{n}}\|\Delta_{q}\textbf{e}_{t}\|^2_{L^2}\Big)\Big\}^{\frac{1}{2}}\nonumber
\end{eqnarray} and the constant $\mu_{4}$ depends only on $A,
\gamma, \tau$ and $\bar{n}$. Furthermore, it is easy to show that
$Q$ satisfies
$$C_{3}\Big(\|U\|_{B^{\sigma}_{2,1}}+\|U_{t}\|_{B^{\sigma-1}_{2,1}}\Big)\leqslant Q\leqslant
C_{4}\Big(\|U\|_{B^{\sigma}_{2,1}}+\|U_{t}\|_{B^{\sigma-1}_{2,1}}\Big)$$
for two positive constants $C_{3},C_{4}$. Choosing
$\delta_{1}=\min\{\frac{\mu_{4}}{2C},
\frac{\bar{\psi}}{(\gamma-1)C}\}$, we complete the proof of
Proposition 5.2 with $\mu_{1}=\frac{\mu_{4}}{2}$.\end{proof}

\noindent\textbf{Proof of Proposition 5.1}. From the assumption
$$
\|U(\cdot,0)\|_{B^{\sigma}_{2,1}}\leqslant \delta_{2},$$ by
Proposition 4.1(local existence), we can determine a time $T_{1}>0
\ (T_{1}<T_{0})$ such that
\begin{eqnarray}\|U(\cdot,t)\|_{B^{\sigma}_{2,1}}+\|U_{t}(\cdot,t)\|_{B^{\sigma-1}_{2,1}}\leqslant
C_{2}\|U(\cdot,t)\|_{B^{\sigma}_{2,1}}\leqslant2C_{2}\delta_{2}
\end{eqnarray} \
for all \ $t\in [0,T_{1}]$.\\ \textbf{Claim}: One can choose a
positive constant $\delta_{2}$ satisfying
$\delta_{2}<\min\{\frac{\delta_{1}}{2C_{2}},\
\frac{\delta_{1}}{C_{1}}\}$ to ensure
\begin{eqnarray}\|U(\cdot,t)\|_{B^{\sigma}_{2,1}}+\|U_{t}(\cdot,t)\|_{B^{\sigma-1}_{2,1}}<
\delta_{1}, \ \ \ \mbox{for all} \ \ t\in [0,T_{0}].\end{eqnarray}
Otherwise, we may assume that there exists a time $T_{2}\
 (T_{1}<T_{2}\leqslant T_{0})$ such that (5.20) is satisfied for all
$t\in[0,T_{2})$ and
\begin{eqnarray}\|U(\cdot,T_{2})\|_{B^{\sigma}_{2,1}}+\|U_{t}(\cdot,T_{2})\|_{B^{\sigma-1}_{2,1}}=
\delta_{1},\end{eqnarray} because (5.20) is satisfied as
$t\in[0,T_{1}]$ for such a choice of $\delta_{2}$.
  By Proposition 5.2, for all
$t\in[0,T^{k}](T^{k}\rightarrow T_{2}-)$,
$$\|U(\cdot,t)\|_{B^{\sigma}_{2,1}}+\|U_{t}(\cdot,t)\|_{B^{\sigma-1}_{2,1}}\leqslant
C_{1}\|U(\cdot,0)\|_{B^{\sigma}_{2,1}}\exp(-\mu_{1}t).$$ In
particular,
\begin{eqnarray*}
\|U(\cdot,T^{k})\|_{B^{\sigma}_{2,1}}+\|U_{t}(\cdot,T^{k})\|_{B^{\sigma-1}_{2,1}}&\leqslant&
C_{1}\|U(\cdot,0)\|_{B^{\sigma}_{2,1}}\exp(-\mu_{1}T^{k})\\
&\leqslant& C_{1}\delta_{2}\exp(-\mu_{1}T^{k})< C_{1}\delta_{2}.
\end{eqnarray*}
By the continuity on $t\in[0,T_{0}]$, we get
$$\|U(\cdot,T_{2})\|_{B^{\sigma}_{2,1}}+\|U_{t}(\cdot,T_{2})\|_{B^{\sigma-1}_{2,1}}\leqslant C_{1}\delta_{2}<\delta_{1},$$
which contradicts (5.21). So, (5.20) holds. From Proposition 4.1
and 5.2, using the boot-strap argument, we can prove Proposition
5.1.\ \ \ \ \ \ \ \ \  \hspace{35mm}$\square$

From the imbedding property in Besov space $B^{\sigma}_{2,1}$,
$(m,\textbf{u},\textbf{e})\in \mathcal{C}^{1}([0,\infty)\times
\mathbb{R}^{N})$ solves (3.1)-(3.2). The choice of $\delta_{1}$ is
sufficient to ensure $\frac{\gamma-1}{2}m+\bar{\psi}>0$. According
to Remark 3.1, we know $(n,\textbf{u},\textbf{e})\in
\mathcal{C}^{1}([0,\infty)\times \mathbb{R}^{N})$ is a solution of
(1.1)-(1.3) with $n>0$.\ Furthermore, we may attain
 Theorem 1.2.

Finally, we show the exponential decay of the vorticity.\\

\noindent\textbf{Proof of Theorem 1.3}. When $N=2$ and $3$, the
$\mathrm{curl}$ of the velocity equation in Eq.(1.1) gives
$$\partial_{t}\omega+\frac{1}{\tau}\omega+\textbf{u}\cdot\nabla \omega-\omega\cdot\nabla\textbf{u}=0\
 \ \ \ \ (\mbox{in particular},\  \omega\cdot\nabla\textbf{u}=0 \  \mbox{when}\ \ N=2 ).$$
Then, we may get
\begin{eqnarray}&&\frac{1}{2}\frac{d}{dt}\|\Delta_{q}\omega\|^2_{L^2}+\frac{1}{\tau}\|\Delta_{q}\omega\|^2_{L^2}\nonumber\\&\leqslant&
C_{5}(\|\nabla\textbf{u}\|_{L^{\infty}}\|\Delta_{q}\omega\|_{L^2}+\|\omega\|_{L^{\infty}}\|\Delta_{q}\nabla\textbf{u}\|_{L^2}\nonumber\\&&+c_{q}\|\nabla
\textbf{u}\|_{B^{\sigma-1}_{2,1}}\|\omega\|_{B^{\sigma-1}_{2,1}})\|\Delta_{q}\omega\|_{L^2}.\end{eqnarray}
Dividing (5.22) by $\|\Delta_{q}\omega\|_{L^2}$ and summing it on
$q\geqslant-1\ (q\in \mathbb{Z})$ after multiplying the factor
$2^{q(\sigma-1)}$, from Theorem 1.2, we have
\begin{eqnarray}
&&\frac{1}{2}\frac{d}{dt}\|\omega(\cdot,t)\|_{B^{\sigma-1}_{2,1}}+\frac{1}{\tau}\|\omega(\cdot,t)\|_{B^{\sigma-1}_{2,1}}\nonumber\\&\leqslant&
C_{5}\|\textbf{u}(\cdot,t)\|_{B^{\sigma}_{2,1}}\|\omega(\cdot,t)\|_{B^{\sigma-1}_{2,1}}\nonumber\\&\leqslant&
C_{5}C_{0}(\|(n_{0}-\bar{n},\textbf{u}_{0},\textbf{e}_{0})\|_{B^{\sigma}_{2,1}})\|\omega(\cdot,t)\|_{B^{\sigma-1}_{2,1}}
\nonumber\\&\leqslant&C_{5}C_{0}\min\Big\{\delta_{0},
\frac{1}{2C_{5}C_{0}\tau}\Big\}\|\omega(\cdot,t)\|_{B^{\sigma-1}_{2,1}}\nonumber\\&\leqslant&\frac{1}{2\tau}\|\omega(\cdot,t)\|_{B^{\sigma-1}_{2,1}}.\end{eqnarray}
Therefore, we obtain the exponential decay of
$\|\omega(\cdot,t)\|_{B^{\sigma-1}_{2,1}}$.\ \ \
\hspace{28mm}$\square$

\section*{Acknowledgment}
This work is supported by NSFC 10571158.

\end{document}